# EARLIEST-DEADLINE-FIRST SERVICE IN HEAVY-TRAFFIC ACYCLIC NETWORKS

By Łukasz Kruk[1], John Lehoczky[2], Steven Shreve[3] and Shu-Ngai Yeung[4]

*Maria Curie-Sklodowska University, Carnegie Mellon University, Carnegie Mellon University and AT&T Laboratories*

This paper presents a heavy traffic analysis of the behavior of multi-class acyclic queueing networks in which the customers have deadlines. We assume the queueing system consists of $J$ stations, and there are $K$ different customer classes. Customers from each class arrive to the network according to independent renewal processes. The customers from each class are assigned a random deadline drawn from a deadline distribution associated with that class and they move from station to station according to a fixed acyclic route. The customers at a given node are processed according to the earliest-deadline-first (EDF) queue discipline. At any time, the customers of each type at each node have a lead time, the time until their deadline lapses. We model these lead times as a random counting measure on the real line. Under heavy traffic conditions and suitable scaling, it is proved that the measure-valued lead-time process converges to a deterministic function of the workload process. A two-station example is worked out in detail, and simulation results are presented to illustrate the predictive value of the theory. This work is a generalization of Doytchinov, Lehoczky and Shreve [*Ann. Appl. Probab.* **11** (2001) 332–379], which developed these results for the single queue case.

**1. Introduction.** The last decade has brought dramatic developments in communication technology. There are now a wide range of types of communication services available and an ever increasing demand for those services. An important component of this demand is for real-time applications,

Received August 2002; revised June 2003.
[1]Supported by the Center for Nonlinear Analysis (NSF Grant DMS-98-03791).
[2]Supported by ONR N00014-92-J-1524 and DARPA F33615-00-C-1729.
[3]Supported by NSF Grant DMS-98-02464 and DMS-01-03814.
[4]Supported by ONR N00014-92-J-1524.

that is, applications with specific timing requirements. Examples are video-conferencing and video-on-demand, in which the timely delivery of packets must be maintained to ensure continuity of the image and sound. Networks servicing real-time applications also arise in production systems where the orders have due dates or in control systems where there are upper bounds on the latency between the occurrence of an event and the control system's response to it.

For queueing networks, the measures of performance are often related to system stability, to queue lengths (e.g., the adequacy of buffer space or the amount of work-in-process in those buffers) or to waiting times (e.g., the delay at a node in the system). For real-time applications, in addition to standard measures of stability, queue length and delay, one must also be concerned whether individual applications are meeting their timing requirements. Determining the waiting time distribution is not sufficient to determine whether a particular queue scheduling policy can satisfy real-time application (customer) timing requirements under various workload conditions.

To model customer timing requirements, we assume that each customer arriving to the system has an initial lead time $\ell$. If the customer arrives at time $t$, then its deadline is met if it exits the network not later than $t + \ell$. To determine whether customers meet their timing requirements, one must dynamically keep track of each customer's lead time, where the lead time is the time remaining until the deadline elapses, that is,

$$\text{lead time } = \text{ initial lead time } - \text{ time elapsed since arrival.}$$

In the study of real-time systems three different types of timing requirements are usually distinguished: hard, firm and soft deadlines. Hard deadlines must be met or a system failure is considered to occur. These applications arise in many control systems, especially in avionics systems or automobile engine control applications. For a computer system to meet hard deadlines, there must be essentially no stochastic aspects associated with the arrival or servicing of an application or these quantities must be bounded above. The subsequent analysis is based on those bounds with service times taking their longest possible value and interarrival times taking their shortest possible value. This worst case approach can result in systems functioning at very low levels of average case utilization to ensure they can meet application deadlines under worst case conditions.

Applications with firm or soft deadlines are permitted to miss their deadlines, usually with some bound on the size or rate of misses. This class of examples include audio and video transmissions, where the end user is able to tolerate a small lack of continuity in the sound or image being transmitted. A firm deadline is one which can be missed and there is no value in completing a task whose deadline has expired, while a soft deadline permits lateness and uses the task completed after the deadline. In this paper we study the soft deadline case but wish to control the fraction of customers that will miss their deadlines and to model the amount of lateness as a function of the workload and the scheduling policies used.

As noted earlier, to study queueing systems in which the customers have deadlines, one must attach a lead-time variable to each customer in the system. It is convenient to model the vector of customer lead times at any time $t$ as a counting measure on $\mathbb{R}$ with a unit atom at the current lead time of each customer and total mass equal to the number of customers in the system at time $t$. Exact analysis of such a measure-valued process is intractable; however, a heavy traffic analysis can be done. Doytchinov, Lehoczky and Shreve (DLS) [5] studied the single queue case in which the customers are processed according to the earliest-deadline-first (EDF) or the first-in-first-out (FIFO) queue disciplines. DLS proved that under heavy traffic conditions, a suitably scaled version of the random lead-time measure converges to a nonrandom function of the limit of the scaled workload process, which in the case of a $GI/G/1$ queue, is a reflected Brownian motion with drift. This paper is focused on generalizing the results of DLS to acyclic queueing networks.

There is some other recent work on heavy-traffic approximations for systems that handle real-time applications. Van Mieghem [15] studied a single server multiclass queueing system. He considered control policies to minimize the total delay cost incurred by customers over a finite time horizon. Markowitz and Wein [12] studied the single machine scheduling problem in a manufacturing context that included customer due dates along with other model features. Lehoczky [10] gave an informal analysis of a single $M/M/1$ queue by constructing the generator for the lead-time vector and showing a scaled version converged to a deterministic limit under the EDF or the processor sharing queue discipline. Lehoczky [9] used these results to study the performance of a threshold access control policy to reduce customer lateness. Lehoczky [11] informally extended the analysis to Jackson networks.

In this paper we extend the approach and the results of DLS to the case of acyclic queueing networks. We assume that a queueing system consists of $J$ stations, and there are $K$ different customer classes. Customers from each class arrive to the queueing network according to independent renewal processes. The customers from each class are assigned a random deadline drawn from some deadline distribution associated with that class, then each moves from station to station according to a fixed route. The $K$ different routes are acyclic, meaning that they visit any of the $J$ stations at most once. Upon completion of their route, they exit the system. If the lead time of a customer at the time of exit is negative, then that customer is late. While lateness is permitted, we wish to determine (and ultimately to control) the fraction of customers that will exit late. Each station will process some subset of the $K$ different customer classes. At each station, customers are queued in lead-time order and preemption (preempt-resume) is allowed. We assume there is no overhead associated with preemption. In the special case that all customers are assigned zero initial lead time, the EDF policy used in this paper becomes first-in-system-first-out (FISFO). Our analysis thus provides information about the time customers spend at various stations and in the system when FISFO is used.

We study the behavior of these acyclic networks under heavy traffic conditions. That is, we consider a sequence of queueing systems in which the traffic intensities at each node approach 1. We prove that if the suitably scaled $J$-dimensional workload process converges appropriately (see Assumption 2.1 and following discussion), then a suitably scaled version of the $K$-dimensional lead-time measure process converges to a deterministic function of the workload.

This paper is organized as follows. Section 2 presents the model, notation and assumptions. Section 3 presents the measure-valued processes associated with customer lead times and the concept of the frontier processes. This section also states the relationship between the limiting workloads and the limiting frontiers, a major result of this paper, which is proved in Section 6. The frontier process was defined in DLS for the single-queue case and formed the basic methodology used to analyze EDF and FIFO queues. This methodology is generalized to the network case in Section 4, which provides several technical results needed for the analysis. Section 5 shows how the equations which provide limiting workloads as a function of limiting frontiers can be inverted, so that one can determine the frontiers as a function of the workloads. Section 7 presents a simple but interesting worked-out example. It also presents simulations which illustrate the accuracy of the heavy traffic theory.

## 2. The model.

2.1. *System topology.* We consider a sequence of queueing systems indexed by $n$. It is assumed that each system is composed of $J$ stations, indexed by 1 through $J$, and $K$ customer classes, indexed by 1 through $K$. Each customer class has a fixed route through the network of stations. Customers in class $k$, $k = 1, \ldots, K$, arrive to the system according to a renewal process, independently of the arrivals of the other customer classes. These customers move through the network, never visiting a station more than once, until they eventually exit the system. However, different customer classes may visit stations in different orders; the system is not necessarily "feed-forward." We define the *path of class $k$ customers* as the sequence of servers they encounter along their way through the network and denote it by

$$(2.1) \qquad \mathcal{P}(k) = (j_{k,1}, j_{k,2}, \ldots, j_{k,m(k)}).$$

In particular, class $k$ customers enter the system at station $j_{k,1}$ and leave it through station $j_{k,m(k)}$. If $j$ is a member of the list of station indices in $\mathcal{P}(k)$, we shall write $j \in \mathcal{P}(k)$.

For $j = 1, \ldots, J$, we define

$$(2.2) \qquad \mathcal{C}(j) \triangleq \{\text{Indices of customer classes that visit station } j\},$$

$$(2.3) \qquad \mathcal{K}_0(j) \triangleq \left\{ \begin{array}{l} \text{Indices of customer classes that enter} \\ \text{station } j \text{ from outside the system} \end{array} \right\}.$$

We assume that for every $j \in \{1, \ldots, J\}$, $\mathcal{C}(j) \neq \varnothing$. For $j = 1, \ldots, J$ and $k \in \mathcal{C}(j)$, we further define

$$(2.4) \qquad \mathcal{S}(k|j) \triangleq \left\{ \begin{array}{l} \text{Indices of stations visited by customer class} \\ k \text{ before visiting station } j, \text{ or } \varnothing \text{ if customer} \\ \text{class } k \text{ enters the system at station } j \end{array} \right\}.$$

We assume that the network is connected, in the sense that for any two stations, there is a way of reaching one station from the other by following fragments of paths of the form (2.1), not necessarily in the forward direction. The network topology captured by (2.1)–(2.4) does not depend on $n$.

2.2. *Arrival times, service times and lead times.* The *customer inter-arrival times* are a sequence of strictly positive, independent and identically distributed random variables $u_k^{i,(n)}$, $i = 1, 2, \ldots$, where the subscript $k$ indicates the customer class and the superscript $(n)$ indexes the queueing system. These are independent across $k$ as well as $i$. We assume that

$$(2.5) \qquad \lambda_k^{(n)} \triangleq (\mathbb{E}\, u_k^{i,(n)})^{-1}, \qquad \alpha_k^{(n)} \triangleq (\operatorname{Var} u_k^{i,(n)})^{1/2}$$

are both defined and finite.

The *customer service times* are a sequence of strictly positive, independent and identically distributed random variables $v_{0,k,j}^{i,(n)}$, $i = 1, 2, \ldots$, where $j$ indicates the station of service and $k \in \mathcal{C}(j)$ indicates the class of customer. The superscript $i$ indexes the order of arrival of customers of class $k$ to the system, which may be different from the order of arrival of class $k$ customers to station $j$. The random variables $v_{0,k,j}^{i,(n)}$ are independent across $k$ and $j$ as well as $i$, and they are independent of the inter-arrival times. We denote by $v_{k,j}^{i,(n)}$ the service times of customers of class $k$ at station $j$, with index $i$ indicating the order of arrival to station $j$. The random variables $\{v_{k,j}^{i,(n)}\}_{i=1}^{\infty}$ are thus a random permutation of the random variables $\{v_{0,k,j}^{i,(n)}\}_{i=1}^{\infty}$. Under the EDF service discipline described below, the index $i$ of arrival at station $j$ of a customer of class $k$ is independent of $v_{k,j}^{i,(n)}$. Therefore, $v_{k,j}^{i,(n)}$, $i = 1, 2, \ldots$, are also independent and identically distributed, with the same distribution as the random variables $v_{0,k,j}^{i,(n)}$. For $j = 1, \ldots, J$ and $k \in \mathcal{C}(j)$, we assume that

$$(2.6) \qquad \begin{aligned} \mu_{k,j}^{(n)} &\triangleq (\mathbb{E} v_{0,k,j}^{i,(n)})^{-1} = (\mathbb{E} v_{k,j}^{i,(n)})^{-1}, \\ \beta_{k,j}^{(n)} &\triangleq (\operatorname{Var} v_{0,k,j}^{i,(n)})^{1/2} = (\operatorname{Var} v_{k,j}^{i,(n)})^{1/2} \end{aligned}$$

are both defined and finite.

Each customer in class $k$ arrives with an *initial lead time* $L_k^{i,(n)}$ having distribution

$$(2.7) \qquad \mathbb{P}\{L_k^{i,(n)} \leq \sqrt{n}y\} = G_k(y),$$

where $G_k$ is a cumulative distribution function. These lead-times are independent and identically distributed across $i$. They are independent across $k$ and independent of the interarrival and service times. Note that $G_k$ does not depend on $n$; the lead-time distributions dilate at rate $\sqrt{n}$ as $n \to \infty$. We assume that for $k = 1, \ldots, K$,

$$(2.8) \qquad y_k^* \triangleq \sup\{y \in \mathbb{R}; G_k(y) < 1\} < \infty.$$

We also assume that for every $n$, the sequences $\{u_k^{i,(n)}\}$, $\{v_{0,k}^{i,(n)}\}$ and $\{L_k^{i,(n)}\}$ are mutually independent over $j \in \{1, \ldots, J\}$, $k \in \mathcal{C}(j)$ and $i = 1, 2, \ldots$.

2.3. *EDF discipline.* Customers are served at each station according to the EDF discipline. That is, the customer with the shortest remaining lead time, regardless of class, is selected for service. We give the analysis for the case of no preemption. If preemption is permitted and we assume preempt-resume, then an obvious simplification of the analysis gives the same results. There is no set-up, switch-over or other type of overhead. Late customers (those with negative lead times) stay in the system until served to completion. We assume the system is empty at time zero.

2.4. *Unscaled queueing processes.* For each station $j = 1, \ldots, J$ and customer class $k = 1, \ldots, K$, we define

$$(2.9) \quad S_k^{m,(n)} \triangleq \sum_{i=1}^{m} u_k^{i,(n)}$$

$$= \text{Time of arrival to the system of the } m\text{th customer of class } k,$$

$$(2.10) \quad A_k^{(n)}(t) \triangleq \max\{m; S_k^{m,(n)} \leq t\}$$

$$= \text{Number of class } k \text{ arrivals to the system by time } t,$$

(2.11) $A_{k,j}^{(n)}(t) \triangleq$ Number of class $k$ arrivals to station $j$ by time $t$,

(2.12) $Q_{k,j}^{(n)}(t) \triangleq$ Number of class $k$ customers at station $j$ at time $t$,

(2.13) $Q_j^{(n)}(t) \triangleq \sum_{k \in \mathcal{C}(j)} Q_{k,j}^{(n)}(t)$

$\qquad$ = Number of customers at station $j$ at time $t$,

(2.14) $V_{0,k,j}^{(n)}(t) \triangleq \sum_{i=1}^{\lfloor t \rfloor} v_{0,k,j}^{i,(n)}$

$\qquad$ = Work for station $j$ associated with the first $\lfloor t \rfloor$ customers of class $k$ to arrive to the system,

(2.15) $V_{k,j}^{(n)}(t) \triangleq \sum_{i=1}^{\lfloor t \rfloor} v_{k,j}^{i,(n)}$

$\qquad$ = Work associated with the first $\lfloor t \rfloor$ customers of class $k$ to arrive at station $j$.

(Here and elsewhere we use the convention $\sum_{i=1}^{0} = \sum_{k \in \varnothing} = 0$.) We further define the *netput at station $j$* to be

(2.16) $$N_j^{(n)}(t) \triangleq \sum_{k \in \mathcal{C}(j)} V_{k,j}^{(n)}(A_{k,j}^{(n)}(t)) - t,$$

and the *cumulative idleness at station $j$* to be $I_j^{(n)}(t) \triangleq -\inf_{0 \le s \le t} N_j^{(n)}(s)$, which is nonnegative because $N_j(0) = 0$. Finally, the *workload at station $j$* is the amount of time it would take to serve all customers at station $j$ to completion, assuming no new customers arrive, and this is $W_j^{(n)}(t) \triangleq N_j^{(n)}(t) + I_j^{(n)}(t)$.

2.5. *Scaled queueing processes.* In order to obtain a limit as $n \to \infty$, it is necessary to scale and sometimes center the above processes. We define

$\widehat{A}_k^{(n)}(t) \triangleq \dfrac{1}{\sqrt{n}}[A_k^{(n)}(nt) - \lambda_k^{(n)} nt], \qquad \widehat{A}_{k,j}^{(n)}(t) \triangleq \dfrac{1}{\sqrt{n}}[A_{k,j}^{(n)}(nt) - \lambda_k^{(n)} nt],$

$\widehat{Q}_{k,j}^{(n)}(t) \triangleq \dfrac{1}{\sqrt{n}} Q_{k,j}^{(n)}(nt), \qquad \widehat{Q}_j^{(n)}(t) \triangleq \sum_{k \in \mathcal{C}(j)} \widehat{Q}_{k,j}^{(n)}(t),$

$\widehat{V}_{0,k,j}^{(n)}(t) \triangleq \dfrac{1}{\sqrt{n}} \sum_{i=1}^{\lfloor nt \rfloor} \left( v_{0,k,j}^{i,(n)} - \dfrac{1}{\mu_{k,j}^{(n)}} \right), \qquad \widehat{V}_{k,j}^{(n)}(t) \triangleq \dfrac{1}{\sqrt{n}} \sum_{i=1}^{\lfloor nt \rfloor} \left( v_{k,j}^{i,(n)} - \dfrac{1}{\mu_{k,j}^{(n)}} \right),$

$\widehat{W}_j^{(n)}(t) \triangleq \dfrac{1}{\sqrt{n}} W_j^{(n)}(nt).$

2.6. *Heavy traffic and convergence assumptions.* For $j = 1, \ldots, J$ and $k \in \mathcal{C}(j)$, the *traffic intensity of class $k$ customers at station $j$* is $\rho_{k,j}^{(n)} \triangleq \lambda_k^{(n)}/\mu_{k,j}^{(n)}$, and the traffic intensity at station $j$ is $\rho_j^{(n)} \triangleq \sum_{k \in \mathcal{C}(j)} \rho_{k,j}^{(n)}$. It is assumed that for all $j$,

$$(2.17) \qquad \gamma_j \triangleq \lim_{n \to \infty} \sqrt{n}(1 - \rho_j^{(n)})$$

exists. Furthermore, it is assumed that for all $k$ and $j$ satisfying $k \in \mathcal{C}(j)$,

$$(2.18) \qquad \begin{aligned} \lambda_k &\triangleq \lim_{n \to \infty} \lambda_k^{(n)}, & \mu_{k,j} &\triangleq \lim_{n \to \infty} \mu_{k,j}^{(n)}, \\ \alpha_k &\triangleq \lim_{n \to \infty} \alpha_k^{(n)}, & \beta_{k,j} &\triangleq \lim_{n \to \infty} \beta_{k,j}^{(n)} \end{aligned}$$

are all defined and $\lambda_k$ and $\mu_{k,j}$ are positive. We define the limiting traffic intensities $\rho_{k,j} \triangleq \lambda_k/\mu_{k,j}$ and $\rho_j \triangleq \sum_{k \in \mathcal{C}(j)} \rho_{k,j}$. We impose the usual Lindeberg condition on the inter-arrival and service times: for $j = 1, \ldots, J$ and $k \in \mathcal{C}(j)$,

$$(2.19) \qquad \begin{aligned} &\lim_{n \to \infty} \mathbb{E}\Big[(u_k^{i,(n)} - (\lambda_k^{(n)})^{-1})^2 \mathbb{I}_{\{|u_k^{i,(n)} - (\lambda_k^{(n)})^{-1}| > c\sqrt{n}\}}\Big] \\ &= \lim_{n \to \infty} \mathbb{E}\Big[(v_{k,j}^{i,(n)} - (\mu_{k,j}^{(n)})^{-1})^2 \mathbb{I}_{\{|v_{k,j}^{i,(n)} - (\mu_{k,j}^{(n)})^{-1}| > c\sqrt{n}\}}\Big] \\ &= 0 \qquad \forall c > 0. \end{aligned}$$

In what follows, the symbol $\Rightarrow$ denotes weak convergence of measures on the space $D_S[0, \infty)$ of right-continuous functions with left limits from $[0, \infty)$ to a Polish space $S$. The topology on this space is a generalization of the topology introduced by Skorokhod for $D_S[0, 1]$. See [2] for details. We take $S = \mathbb{R}$ (or $\mathbb{R}^d$, with appropriate dimension $d$, for vector-valued functions) unless explicitly stated otherwise.

Theorem 3.1 [14], together with (2.19) and the independence assumptions of Section 2.2, implies that for every $j = 1, \ldots, J$ and $k \in \mathcal{C}(j)$ and every $y \leq y_k^*$, we have

$$(2.20) \quad \widehat{T}_{0,k,j}^{(n)}(t; y) \triangleq \frac{1}{\sqrt{n}} \sum_{i=1}^{\lfloor nt \rfloor} \Big[ v_{0,k,j}^{i,(n)} \mathbb{I}_{\{L_k^{i,(n)} \leq \sqrt{n}\, y\}} - \frac{1}{\mu_{k,j}^{(n)}} G_k(y) \Big] \Rightarrow T_{k,j}^*(t; y),$$

where $T_{k,j}^*(t; y)$ is continuous in $t$. Putting $y = y_k^*$ into (2.20) and using the fact that the sequences $\{v_{0,k,j}^{i,(n)}\}_{i=1}^{\infty}$ and $\{v_{k,j}^{i,(n)}\}_{i=1}^{\infty}$ have the same distribution, we conclude that for $j = 1, \ldots, J$ and $k \in \mathcal{C}(j)$,

$$(2.21) \qquad \widehat{V}_{0,k,j}^{(n)} \Rightarrow \widehat{V}_{k,j}^*, \qquad \widehat{V}_{k,j}^{(n)} \Rightarrow \widehat{V}_{k,j}^*,$$

where $\widehat{V}_{k,j}^*$ is a continuous process. In fact, if $\beta_{k,j} > 0$, then $\widehat{V}_{k,j}^*$ is a Brownian motion. Similarly, Theorem 3.1 of [14] and Theorem 14.6 of [2] imply that, for every $k$, there exists a continuous process $A_k^*$ such that

$$(2.22) \qquad \widehat{A}_k^{(n)} \Rightarrow A_k^*.$$

By (2.21), (2.22) and a standard argument (see, e.g., [8], Corollary 3.2), we have

$$(2.23) \qquad \widehat{M}_{0,k,j}^{(n)}(t) \triangleq \frac{1}{\sqrt{n}}[V_{0,k,j}^{(n)}(A_k^{(n)}(nt)) - n\rho_{k,j}^{(n)}t] \Rightarrow M_{k,j}^*(t),$$

where $M_{k,j}^*$ is continuous. We also make the following convergence assumption:

ASSUMPTION 2.1. For every $j$ and $k \in \mathcal{C}(j)$, there exists a continuous process $A_{k,j}^*$ such that

$$(2.24) \qquad \widehat{A}_{k,j}^{(n)} \Rightarrow A_{k,j}^*.$$

There exists a $J$-dimensional continuous process $(W_1^*, W_2^*, \ldots, W_J^*)$ such that

$$(2.25) \qquad (\widehat{W}_1^{(n)}, \widehat{W}_2^{(n)}, \ldots, \widehat{W}_J^{(n)}) \Rightarrow (W_1^*, W_2^*, \ldots, W_J^*).$$

In feed-forward networks, (2.24) and (2.25) hold under FIFO (see [13]) and EDF (see [17]). Because our network is not of the feed-forward type, there are no known general conditions which guarantee (2.24) and (2.25). However, the literature contains a number of special cases of our model in which (2.24) and (2.25) hold. Rather than take one of these special cases as a starting point, we choose to begin with Assumption 2.1 because this is all we shall need in order to obtain convergence of scaled lead-time profiles.

## 3. Measure-valued processes and frontiers.
Our goal is to obtain a characterization of the lead-time profiles of the customers queued at the $J$ stations in the system in terms of the limiting workload process $(W_1^*, W_2^*, \ldots, W_J^*)$ in (2.25). These lead-time profile processes are measure-valued. More precisely, they take values in the space $\mathcal{M}$ of finite, nonnegative measures on $\mathcal{B}(\mathbb{R})$, the Borel $\sigma$-algebra on $\mathbb{R}$, equipped with the weak topology. In what follows, we shall denote by $\mathcal{M}^J$ the $J$-fold product of $\mathcal{M}$ (with the product topology). For a Borel set $B \subset \mathbb{R}$, we set

$$(3.1) \qquad \mathcal{W}_{k,\ell}^{j,(n)}(t)(B) \triangleq \left\{ \begin{array}{l} \text{Work for station } \ell \text{ represented by} \\ \text{class } k \text{ customers at station } j \text{ with} \\ \text{lead times in } B \text{ at time } t \end{array} \right\},$$

$$(3.2) \qquad \mathcal{W}_\ell^{j,(n)}(t)(B) \triangleq \left\{ \begin{array}{l} \text{Work for station } \ell \text{ represented by customers} \\ \text{at station } j \text{ with lead times in } B \text{ at time } t \end{array} \right\}.$$

Then $\mathcal{W}_{k,j}^{(n)}(t)(B) \triangleq \mathcal{W}_{k,j}^{j,(n)}(t)(B)$ is the work at station $j$ represented by class $k$ customers at that station with lead times in $B$ at time $t$, and $\mathcal{W}_j^{(n)}(t)(B) \triangleq \mathcal{W}_j^{j,(n)}(t)(B)$ is the work at station $j$ corresponding to all customers at that station with lead times in $B$ at time $t$. We also define

$$(3.3) \qquad \mathcal{Q}_j^{(n)}(t)(B) \triangleq \left\{ \begin{array}{l} \text{Number of customers at station } j \\ \text{with lead times in } B \text{ at time } t \end{array} \right\},$$

$$(3.4) \qquad \mathcal{A}_k^{(n)}(t)(B) \triangleq \left\{ \begin{array}{l} \text{Number of class } k \text{ customers arriving} \\ \text{to the system by time } t \text{ and having} \\ \text{lead times at time } t \text{ in } B, \text{ whether or} \\ \text{not still in the system at time } t \end{array} \right\}$$

and

$$(3.5) \, \mathcal{V}_{0,k,j}^{(n)}(t)(B) \triangleq \left\{ \begin{array}{l} \text{Work for station } j \text{ associated with customers} \\ \text{of type } k \text{ arriving to the system by time } t \text{ and} \\ \text{having lead times at time } t \text{ in } B, \text{ whether or} \\ \text{not still in the system at time } t \end{array} \right\}.$$

The scaled versions of these processes are

$$\widehat{\mathcal{W}}_{k,\ell}^{j,(n)}(t)(B) \triangleq \frac{1}{\sqrt{n}} \mathcal{W}_{k,\ell}^{j,(n)}(nt)(\sqrt{n}B),$$

$$\widehat{\mathcal{W}}_{\ell}^{j,(n)}(t)(B) \triangleq \frac{1}{\sqrt{n}} \mathcal{W}_{\ell}^{j,(n)}(nt)(\sqrt{n}B),$$

$$\widehat{\mathcal{W}}_{k,j}^{(n)}(t)(B) \triangleq \frac{1}{\sqrt{n}} \mathcal{W}_{k,j}^{(n)}(nt)(\sqrt{n}B),$$

$$\widehat{\mathcal{W}}_j^{(n)}(t)(B) \triangleq \frac{1}{\sqrt{n}} \mathcal{W}_j^{(n)}(nt)(\sqrt{n}B),$$

$$\widehat{\mathcal{Q}}_j^{(n)}(t)(B) \triangleq \frac{1}{\sqrt{n}} \mathcal{Q}_j^{(n)}(nt)(\sqrt{n}B),$$

$$\widehat{\mathcal{A}}_k^{(n)}(t)(B) \triangleq \frac{1}{\sqrt{n}} \mathcal{A}_k^{(n)}(nt)(\sqrt{n}B),$$

$$\widehat{\mathcal{V}}_{0,k,j}^{(n)}(t)(B) \triangleq \frac{1}{\sqrt{n}} \mathcal{V}_{0,k,j}^{(n)}(nt)(\sqrt{n}B).$$

We introduce *frontier processes*

$$(3.6) \quad F_{k,j}^{(n)}(t) \triangleq \left\{ \begin{array}{l} \text{Largest lead time of any class } k \text{ customer} \\ \text{who has ever been in service at station } j, \\ \text{or } \sqrt{n}y_k^* - t \text{ if no such customer exists or} \\ \text{if this quantity is larger than the former one} \end{array} \right\},$$

$$(3.7) \quad F_j^{(n)}(t) \triangleq \max_{k \in \mathcal{C}(j)} F_{k,j}^{(n)}(t).$$

The scaled versions of these processes are

$$(3.8) \qquad \widehat{F}_{k,j}^{(n)}(t) = \frac{1}{\sqrt{n}} F_{k,j}^{(n)}(nt), \qquad \widehat{F}_{j}^{(n)}(t) \triangleq \frac{1}{\sqrt{n}} F_{j}^{(n)}(nt).$$

The next step is to define a set $D$ which contains the $J$-dimensional vector-valued process $(\widehat{F}_1^{(n)}, \widehat{F}_2^{(n)}, \ldots, \widehat{F}_J^{(n)})$. To do this, we begin with a permutation $\pi = (\pi_1, \pi_2, \ldots, \pi_J)$ of the integers $(1, 2, \ldots, J)$. Given such a permutation and an integer $m \in \{1, \ldots, J\}$, we define

$$(3.9) \qquad \mathcal{K}_{m-1}^{\pi}(j) \triangleq \{k \in \mathcal{C}(j); \mathcal{S}(k|j) \subset \{\pi_1, \ldots, \pi_{m-1}\}\},$$

$$(3.10) \qquad \mathcal{J}_{m-1}^{\pi} \triangleq \{j; \mathcal{K}_{m-1}^{\pi}(j) \neq \varnothing\} \setminus \{\pi_1, \ldots, \pi_{m-1}\}.$$

By convention, if $m = 1$, then $\{\pi_1, \ldots, \pi_{m-1}\} = \varnothing$, $\mathcal{K}_0^{\pi}(j) = \mathcal{K}_0(j)$, the set of indices of customer classes which enter the system at station $j$, and $\mathcal{J}_0^{\pi} = \mathcal{J}_0$, the set of stations which serve as the entry point for at least one external arrival process. These two sets do not depend on the permutation $\pi$. Subsequent sets do. The set $\mathcal{K}_{m-1}^{\pi}(j)$ is the set of all customer classes that visit station $j$ and visit only stations in the set $\{\pi_1, \ldots, \pi_{m-1}\}$ before arriving at station $j$. We say it is the *set of customer classes which reach station $j$ through $\{\pi_1, \ldots, \pi_{m-1}\}$*. The set $\mathcal{J}_{m-1}^{\pi}$ is the set of all stations $j$ not in the set $\{\pi_1, \ldots, \pi_{m-1}\}$ which are visited by at least one customer class of the type just described. We say that $\mathcal{J}_{m-1}^{\pi}$ is the *set of stations which can be reached through $\{\pi_1, \ldots, \pi_{m-1}\}$*. Note that both $\mathcal{K}_{m-1}^{\pi}(j)$ and $\mathcal{J}_{m-1}^{\pi}$ depend only on $(\pi_1, \ldots, \pi_{m-1})$, not the full permutation. Thus, we shall sometimes write $\mathcal{K}_{m-1}^{(\pi_1, \ldots, \pi_{m-1})}(j)$ and $\mathcal{J}_{m-1}^{(\pi_1, \ldots, \pi_{m-1})}$ instead of $\mathcal{K}_{m-1}^{\pi}(j)$ and $\mathcal{J}_{m-1}^{\pi}$. Finally, we define

$$(3.11) \qquad \begin{aligned} \Pi \triangleq \{&\pi; \pi \text{ is a permutation of } 1, \ldots, J \\ &\text{and } \pi_m \in \mathcal{J}_{m-1}^{\pi} \text{ for all } m = 1, \ldots, J\}. \end{aligned}$$

In other words, $\Pi$ is the set of all permutations $\pi = (\pi_1, \ldots, \pi_J)$ such that, for each $m$, the station $\pi_m$ can be reached through $\{\pi_1, \ldots, \pi_{m-1}\}$. For $\pi \in \Pi$, we set

$$(3.12) \quad D^{\pi} \triangleq \left\{ y \in \mathbb{R}^J; y_{\pi_1} \geq \cdots \geq y_{\pi_J} \text{ and } y_{\pi_m} \leq \max_{k \in \mathcal{K}_{m-1}^{\pi}(\pi_m)} y_k^* \ \forall m \right\},$$

$$(3.13) \quad D \triangleq \bigcup_{\pi \in \Pi} D^{\pi}.$$

LEMMA 3.1. *For all $t \geq 0$, the random vector $(\widehat{F}_1^{(n)}(t), \widehat{F}_2^{(n)}(t), \ldots, \widehat{F}_J^{(n)}(t))$ takes values in the set $D$.*

PROOF. We must construct a permutation $\pi = (\pi_1, \ldots, \pi_J) \in \Pi$ such that $\widehat{F}_{\pi_1}^{(n)}(t) \geq \widehat{F}_{\pi_2}^{(n)}(t) \geq \cdots \geq \widehat{F}_{\pi_J}^{(n)}(t)$ and $\widehat{F}_{\pi_m}^{(n)}(t) \leq \max_{k \in \mathcal{K}_{m-1}^{\pi}(\pi_m)} y_k^*$ for every $m$. We do this by induction.

We note first that because customer class $k$ visits consecutive stations $j_{k,1}, j_{k,2}, \ldots, j_{k,m(k)}$ in $\mathcal{P}(k)$,

$$(3.14) \qquad \sqrt{n}y_k^* \geq F_{k,j_{k,1}}^{(n)}(nt) \geq F_{k,j_{k,2}}^{(n)}(nt) \geq \cdots \geq F_{k,j_{k,m(k)}}^{(n)}(nt).$$

This implies that the largest frontier must be at a station which is in $\mathcal{J}_0$, the set of stations that have arrivals from outside the system. We select a station $\pi_1 \in \mathcal{J}_0$ whose frontier $F_{\pi_1}^{(n)}(nt)$ is maximal. If this maximal frontier is the lead time of a customer which has been in service, we may choose $\pi_1$ to be the station where that customer entered the system. If the maximal frontier is not the lead time of a customer who has been in service, then it is of the form $\sqrt{n}y_{\bar{k}}^* - nt$ for some customer class $\bar{k}$. In this case, we choose $\pi_1$ to be the station where this customer class enters the system, so that $\bar{k}$ is in $\mathcal{K}_0(\pi_1)$, the set of customer classes that enter the system at station $\pi_1$. In either case, we obtain

$$(3.15) \qquad \max_{k \in \mathcal{K}_0(\pi_1)} \sqrt{n}y_k^* \geq F_{\pi_1}^{(n)}(nt) \geq \max_{j \neq \pi_1} F_j^{(n)}(nt).$$

For the induction hypothesis, we assume for some $m \in \{2, \ldots, J\}$ that we have constructed $\pi_1, \ldots, \pi_{m-1}$ such that:

(i) for each $i \leq m-1$, station $\pi_i$ is reached through $\{\pi_1, \ldots, \pi_{i-1}\}$,

(ii) for each $i \leq m-1$,

$$(3.16) \qquad \max_{k \in \mathcal{K}_{i-1}^{(\pi_1,\ldots,\pi_{i-1})}(\pi_i)} \sqrt{n}y_k^* \geq F_{\pi_i}^{(n)}(nt),$$

(iii) we have

$$(3.17) \qquad F_{\pi_1}^{(n)}(nt) \geq \cdots \geq F_{\pi_{m-1}}^{(n)}(nt) \geq \max_{j \notin \{\pi_1,\ldots,\pi_{m-1}\}} F_j^{(n)}(nt).$$

If the maximal frontier among $F_j^{(n)}(nt)$ for $j \notin \{\pi_1, \ldots, \pi_{m-1}\}$ is the lead time of a customer that has been in service, we may choose $\pi_m$ to be the station where that customer first reaches a station not in $\{\pi_1, \ldots, \pi_{m-1}\}$. If this maximal frontier is not the lead time of a customer who has been in service, then it is of the form $\sqrt{n}y_{\bar{k}}^* - nt$ for some customer class $\bar{k}$. In this case, we choose $\pi_m$ to be the first station not in $\{\pi_1, \ldots, \pi_{m-1}\}$ reached by this customer class. In either case, we obtain

$$(3.18) \qquad \max_{k \in \mathcal{K}_{m-1}^{(\pi_1,\ldots,\pi_{m-1})}(\pi_m)} \sqrt{n}y_k^* \geq F_{\pi_m}^{(n)}(nt) \geq \max_{j \notin \{\pi_1,\ldots,\pi_m\}} F_j^{(n)}(nt).$$

The induction step is complete.

Once the induction has concluded, we have constructed a permutation satisfying properties (i)–(iii) for $m = J + 1$. Dividing (3.16) and (3.17) by $\sqrt{n}$, we obtain the desired properties for the scaled frontiers.  □

For $k = 1, \ldots, K$, we define

$$(3.19) \qquad H_k(y) \triangleq \int_y^\infty (1 - G_k(x)) \, dx, \qquad y \in \mathbb{R}.$$

This function is strictly decreasing on $(-\infty, y_k^*]$, mapping this half-line onto $[0, \infty)$. We next define $\Phi = (\Phi_1, \ldots, \Phi_J) \colon \mathbb{R}^J \to [0, \infty)^J$ by

$$(3.20) \qquad \Phi_j(y_1, \ldots, y_J) \triangleq \sum_{k \in \mathcal{C}(j)} \rho_{k,j} \Big[ H_k(y_j) - H_k\Big( \min_{i \in \mathcal{S}(k|j)} y_i \Big) \Big]^+,$$
$$j = 1, \ldots, J.$$

In the above definition and in all that follows, the minimum taken over the empty set should be interpreted as $\infty$. The main results of this paper are the following two theorems.

THEOREM 3.2 (Convergence of scaled frontiers). *The function $\Phi$ is a homeomorphism of $D$ onto $[0, \infty)^J$. With*

$$(3.21) \qquad (F_1^*, \ldots, F_J^*) \triangleq \Phi^{-1}(W_1^*, \ldots, W_J^*),$$

*we have*

$$(3.22) \qquad (\widehat{F}_1^{(n)}, \ldots, \widehat{F}_J^{(n)}) \Rightarrow (F_1^*, \ldots, F_J^*).$$

We define $\mathcal{W}^*(t) = (\mathcal{W}_1^*(t), \ldots, \mathcal{W}_J^*(t))$ and $\mathcal{Q}^*(t) = (\mathcal{Q}_1^*(t), \ldots, \mathcal{Q}_J^*(t))$, which take values in $\mathcal{M}^J$, the set of $J$-dimensional vectors of measures on $\mathbb{R}$, by specifying their values on half-lines of the form $(y, \infty)$ for all $y \in \mathbb{R}$. This is done for $j = 1, \ldots, J$ by the formulas

$$(3.23) \quad \mathcal{W}_j^*(t)(y, \infty) \triangleq \sum_{k \in \mathcal{C}(j)} \rho_{k,j} \Big[ H_k(y \vee F_j^*(t)) - H_k\Big( \min_{i \in \mathcal{S}(k|j)} F_i^*(t) \Big) \Big]^+,$$

$$(3.24) \quad \mathcal{Q}_j^*(t)(y, \infty) \triangleq \sum_{k \in \mathcal{C}(j)} \lambda_k \Big[ H_k(y \vee F_j^*(t)) - H_k\Big( \min_{i \in \mathcal{S}(k|j)} F_i^*(t) \Big) \Big]^+.$$

THEOREM 3.3 (Convergence of scaled workloads and queue lengths). *We have*

$$(3.25) \qquad (\widehat{\mathcal{W}}_1^{(n)}, \ldots, \widehat{\mathcal{W}}_J^{(n)}) \Rightarrow \mathcal{W}^*, \qquad (\widehat{\mathcal{Q}}_1^{(n)}, \ldots, \widehat{\mathcal{Q}}_J^{(n)}) \Rightarrow \mathcal{Q}^*.$$

*The weak convergence in* (3.25) *takes place in* $D_{\mathcal{M}^J}[0, \infty)$.

**4. Customers behind the frontiers.**   In this section we prove the crucial observation that both the number of customers at each station $j$ with lead times smaller than or equal to the current frontier $F_j^{(n)}(t)$ and the work for the system associated with these customers are negligible. This is done in several steps, leading to Corollary 4.7. Along the way, we show tightness of the scaled frontier processes (Lemma 4.6). Both these results will be used in Section 6 to prove Theorem 3.2.

PROPOSITION 4.1.   *Let* $j = 1, \dots, J$, $k \in \mathcal{C}(j)$, $-\infty < y_0 < y_k^*$ *and* $T > 0$ *be given. As* $n \to \infty$,

$$(4.1) \qquad \sup_{y_0 \le y \le y_k^*} \sup_{0 \le t \le T} |\widehat{\mathcal{V}}_{0,k,j}^{(n)}(t)(y, \infty) + \rho_{k,j}[H_k(y + \sqrt{n}t) - H_k(y)]| \xrightarrow{P} 0,$$

$$(4.2) \qquad \sup_{y_0 \le y \le y_k^*} \sup_{0 \le t \le T} |\widehat{\mathcal{A}}_k^{(n)}(t)(y, \infty) + \lambda_k[H_k(y + \sqrt{n}t) - H_k(y)]| \xrightarrow{P} 0.$$

The processes in Proposition 4.1 do not take departures into account. Because this proposition is concerned only with arrivals, its proof can be given following the proof of Proposition 3.4 of [5]. We do not repeat that proof here, but instead give a heuristic argument. Let us first consider (4.2), which asserts that asymptotically, the "density" of the measure-valued process $\widehat{\mathcal{A}}_k^{(n)}(t)$ is the same as the density

$$(4.3) \qquad \lambda_k^{(n)}[H_k'(y + \sqrt{n}t) - H_k'(y)] = \lambda_k^{(n)}[G_k(y + \sqrt{n}t) - G_k(y)].$$

In order for a class $k$ customer to have lead time $\hat{y}$ at time $\hat{t}$, it must arrive at some time $\hat{t} - \hat{s}$ prior to $\hat{t}$ and be assigned lead time $\hat{y} + \hat{s}$ upon arrival. If $G_k$ has a density, then the density of the assigned lead-time distribution is $\frac{1}{\sqrt{n}} G_k'\left(\frac{\hat{y}+\hat{s}}{\sqrt{n}}\right)$ [see (2.7)], and multiplying by the arrival rate $\lambda_k^{(n)}$, we obtain the density of class $k$ customers with lead times $\hat{y}$:

$$\frac{\lambda_k^{(n)}}{\sqrt{n}} \int_0^{\hat{t}} G_k'\left(\frac{\hat{y}+\hat{s}}{\sqrt{n}}\right) d\hat{s} = \lambda_k^{(n)}\left[G_k\left(\frac{\hat{y}+\hat{t}}{\sqrt{n}}\right) - G_k\left(\frac{\hat{y}}{\sqrt{n}}\right)\right].$$

The heavy traffic scaling considers the density of $\frac{1}{\sqrt{n}}$ times the actual number of customers whose lead times are $y = \frac{\hat{y}}{\sqrt{n}}$ at scaled times $t = \frac{\hat{t}}{n}$. This density is the right-hand side of (4.3). (The Jacobian $\frac{d\hat{y}}{dy} = \sqrt{n}$ is canceled when we divide the customer count by $\sqrt{n}$.) Under the heavy traffic scaling, the work brought by customers of class $k$ to station $j$ is the average work per customer, $(\mu_{k,j}^{(n)})^{-1}$, times the number of customers. Multiplying the right-hand side of (4.3) by $(\mu_{k,j}^{(n)})^{-1}$, we obtain $\rho_{k,j}^{(n)}[G_k(y + \sqrt{n}t) - G_k(y)]$, which explains (4.1).

COROLLARY 4.2.   *Let* $j = 1, \ldots, J$, $k \in \mathcal{C}(j)$, $-\infty < y_0 < y_k^*$ *and* $T > 0$ *be given. As* $n \to \infty$,

$$(4.4) \qquad \begin{aligned} \sup_{y_0 \leq y \leq y_k^*} \sup_{0 \leq t \leq T} \widehat{\mathcal{V}}_{0,k,j}^{(n)}(t)\{y\} &\xrightarrow{P} 0, \\ \sup_{y_0 \leq y \leq y_k^*} \sup_{0 \leq t \leq T} \widehat{\mathcal{A}}_k^{(n)}(t)\{y\} &\xrightarrow{P} 0. \end{aligned}$$

This corollary is a consequence of the fact that the limiting measures for $\widehat{\mathcal{V}}_{0,k,j}^{(n)}(t)$ and $\widehat{\mathcal{A}}_k^{(n)}(t)$ have densities and, hence, do not charge points. Its proof is similar to the proof of Corollary 3.5 of [5], and we refer the reader there for details.

Lemma 4.3 and Corollary 4.4 generalize Proposition 3.6 in [5] to the case of acyclic networks. We provide the details of these proofs.

LEMMA 4.3.   *For all* $k$ *and* $j \in \mathcal{P}(k)$, *we have*

$$(4.5) \qquad \widehat{\mathcal{Q}}_j^{(n)}(-\infty, \widehat{F}_{k,j}^{(n)}) \Rightarrow 0, \qquad \widehat{\mathcal{W}}_j^{(n)}(-\infty, \widehat{F}_{k,j}^{(n)}) \Rightarrow 0,$$

*and for every station* $\ell$ *satisfying* $j \in \mathcal{S}(k|\ell)$,

$$(4.6) \qquad \widehat{\mathcal{W}}_\ell^{j,(n)}(-\infty, \widehat{F}_{k,j}^{(n)}) \Rightarrow 0.$$

*The weak convergences in* (4.5) *and* (4.6) *take place in* $D_{\mathbb{R}}[0, \infty)$.

PROOF.   We fix $k_0$ and prove (4.5) and (4.6) by induction along $\mathcal{P}(k_0)$. Let $j_0 \in \mathcal{P}(k_0)$ be given and make the induction hypothesis that (4.5) and (4.6) hold for $k = k_0$ and every $j \in \mathcal{S}(k_0|j_0)$. If $j_0$ is the first station in $\mathcal{P}(k_0)$, this hypothesis is vacuous.

From the induction hypothesis (4.6) we have

$$\widehat{\mathcal{W}}_{j_0}^{j,(n)}(-\infty, \widehat{F}_{k_0,j}^{(n)}) \Rightarrow 0$$

for every $j \in \mathcal{S}(k_0|j_0)$. For such a station $j$, we have $\widehat{F}_{k_0,j}^{(n)} \geq \widehat{F}_{k_0,j_0}^{(n)}$ and, hence,

$$(4.7) \qquad \widehat{\mathcal{W}}_{j_0}^{j,(n)}(-\infty, \widehat{F}_{k_0,j_0}^{(n)}) \Rightarrow 0.$$

We now prove (4.5) and (4.6) for $k = k_0$ and $j = j_0$. Toward this end, we define

$$\tau_{k_0,j_0}^{(n)}(t) = \sup\{s \leq t; \mathcal{Q}_{j_0}^{(n)}(ns)(-\infty, F_{k_0,j_0}^{(n)}(ns)) = 0\}$$
$$= \sup\{s \leq t; \widehat{\mathcal{Q}}_{j_0}^{(n)}(s)(-\infty, \widehat{F}_{k_0,j_0}^{(n)}(s)) = 0\},$$

and note that

$$(4.8) \qquad \mathcal{Q}_{j_0}^{(n)}(n\tau_{k_0,j_0}^{(n)}(t)-)(-\infty, F_{k_0,j_0}^{(n)}(n\tau_{k_0,j_0}^{(n)}(t))) = 0.$$

Indeed, from the definition of $\tau_{k_0,j_0}^{(n)}(t)$,

$$(4.9) \qquad \mathcal{Q}_{j_0}^{(n)}(n\tau_{k_0,j_0}^{(n)}(t)-)(-\infty, F_{k_0,j_0}^{(n)}(n\tau_{k_0,j_0}^{(n)}(t)-)) = 0$$

and the only way in which $F_{k_0,j_0}^{(n)}$ can jump up at time $n\tau_{k_0,j_0}^{(n)}(t)$ is that a customer of class $k_0$ with lead time greater than any customer of this class who has ever been in service at station $j_0$ begins to receive service at time $n\tau_{k_0,j_0}^{(n)}(t)$. But then, by the EDF service discipline,

$$\mathcal{Q}_{j_0}^{(n)}(n\tau_{k_0,j_0}^{(n)}(t))(-\infty, F_{k_0,j_0}^{(n)}(n\tau_{k_0,j_0}^{(n)}(t))) = 0$$

and, consequently, the value of $\mathcal{Q}_{j_0}^{(n)}(\cdot)(-\infty, F_{k_0,j_0}^{(n)}(\cdot))$ remains zero at least until the next arrival to server $j_0$ after time $n\tau_{k_0,j_0}^{(n)}(t)$, which contradicts the definition of $\tau_{k_0,j_0}^{(n)}(t)$. Thus, $F_{k_0,j_0}^{(n)}$ cannot jump up at time $n\tau_{k_0,j_0}^{(n)}(t)$, so (4.8) follows from (4.9). For $\tau_{k_0,j_0}^{(n)}(t) \leq s \leq t$, we have

$$(4.10) \qquad F_{k_0,j_0}^{(n)}(ns) = F_{k_0,j_0}^{(n)}(n\tau_{k_0,j_0}^{(n)}(t)) - n(s - \tau_{k_0,j_0}^{(n)}(t)).$$

The first step is to prove

$$(4.11) \qquad t - \tau_{k_0,j_0}^{(n)}(t) \Rightarrow 0,$$

and, subsequently, to upgrade this convergence to

$$(4.12) \qquad \sqrt{n}(t - \tau_{k_0,j_0}^{(n)}(t)) \Rightarrow 0.$$

The convergences in (4.11) and (4.12) are in $D_{\mathbb{R}}[0,\infty)$. The key inequality is

$$(4.13) \quad \begin{aligned} \mathcal{W}_{j_0}^{(n)}(nt)&(-\infty, F_{k_0,j_0}^{(n)}(nt)) \\ &\leq H^{(n)}(nt) + \sum_{j \in \mathcal{S}(k_0|j_0)} J_j^{(n)}(nt) + D^{(n)}(nt) \\ &\quad + \sum_{\substack{k \in \mathcal{C}(j_0) \\ k \neq k_0}} K_k^{(n)}(nt) - n(t - \tau_{k_0,j_0}^{(n)}(t)) + R^{(n)}(nt), \end{aligned}$$

where the terms on the right-hand side are defined below. The first term,

$$H^{(n)}(nt) \triangleq \mathcal{W}_{j_0}^{(n)}(n\tau_{k_0,j_0}^{(n)}(t))(-\infty, F_{k_0,j_0}^{(n)}(n\tau_{j_0,k_0}^{(n)}(t))),$$

accounts for the work arriving to station $j_0$ by time $n\tau_{k_0,j_0}^{(n)}(t)$. A typical summand in the second term,

$$J_j^{(n)}(nt) \triangleq \mathcal{W}_{j_0}^{j,(n)}(n\tau_{k_0,j_0}^{(n)}(t))(-\infty, F_{k_0,j_0}^{(n)}(n\tau_{k_0,j_0}^{(n)}(t))),$$

is the work for station $j_0$ at upstream station $j$ which is ahead of $F_{k_0,j_0}^{(n)}(n\tau_{k_0,j_0}^{(n)}(t))$ at time $n\tau_{k_0,j_0}^{(n)}$ and, hence, has the potential to arrive at station $j_0$ by time $nt$ ahead of $F_{k_0,j_0}^{(n)}(nt)$. If $j_0$ is the first station in $\mathcal{P}(k_0)$, these terms do not appear. The third term,

$$D^{(n)}(nt) = \sum_{i=1}^{\infty} v_{0,k_0,j_0}^{i,(n)} \mathbb{1}_{\{n\tau_{k_0,j_0}^{(n)}(t) < S_{k_0}^{i,(n)} \le nt\}} \mathbb{1}_{\{L_{k_0}^{i,(n)} - (nt - S_{k_0}^{i,(n)}) \le F_{k_0,j_0}^{(n)}(nt)\}},$$

is the work of type $k_0$ arriving to the system during the time interval $(n\tau_{k_0,j_0}^{(n)}, nt]$ with lead time upon arrival that puts it ahead of $F_{k_0,j_0}(nt)$ at time $nt$. A typical summand in the fourth term,

$$K_k^{(n)}(nt) = \sum_{i=1}^{\infty} v_{k,j_0}^{i,(n)} \mathbb{1}_{\{A_{k,j_0}(n\tau_{k_0,j_0}^{(n)}(t)) < i \le A_{k,j_0}^{(n)}(nt)\}},$$

is the work of type $k \ne k_0$ arriving to station $j_0$ during the time interval $(n\tau_{k_0,j_0}^{(n)}, nt]$. The fifth term, $-n(t - \tau_{k_0,j_0}^{(n)}(t))$, is the work accomplished by the server during the time interval $(n\tau_{k_0,j_0}^{(n)}(t), nt]$, some of which may be devoted to a customer already in service at time $n\tau_{k_0,j_0}^{(n)}(t)$ whose lead time is greater than or equal to $F_{k_0,j_0}^{(n)}(n\tau_{k_0,j_0}^{(n)}(t))$ and the remainder of which at each time $s \in (\tau_{k_0,j_0}^{(n)}, t]$ must be devoted to customers with lead times less than $F_{k_0,j_0}^{(n)}(ns)$. The final term,

$$R^{(n)}(nt) \triangleq \max\{v_{k,j_0}^{i,(n)}; 1 \le i \le A_{k,j_0}^{(n)}(nt), k \in \mathcal{C}(j_0)\},$$

is an upper bound on the amount of work that can be devoted to a customer already in service at time $n\tau_{k_0,j_0}^{(n)}(t)$. If there is preemption, this final term does not appear.

We fix $T > 0$ and estimate the terms appearing on the right-hand side of (4.13). The terms $o(\sqrt{n})$, $O(\sqrt{n})$ and $O(n^{-1/2})$ in the following argument depend on $T$ but not on $t \in [0, T]$. For $t \in [0, T]$, (4.8) implies

$$
\begin{aligned}
(4.14) \qquad H^{(n)}(nt) &\le \max_{0 \le s \le T}[W_{j_0}^{(n)}(ns) - W_{j_0}^{(n)}(ns-)] \\
&= \sqrt{n} \max_{0 \le s \le T}[\widehat{W}_{j_0}^{(n)}(s) - \widehat{W}_{j_0}^{(n)}(s-)] = o(\sqrt{n})
\end{aligned}
$$

because $W_{j_0}^*$ in (2.25) is continuous. Furthermore,

$$
\begin{aligned}
(4.15) \qquad J_j^{(n)}(nt) &\le \max_{0 \le s \le T} \mathcal{W}_{j_0}^{j,(n)}(ns)(-\infty, F_{k_0,j_0}^{(n)}(ns)) \\
&= \sqrt{n} \max_{0 \le s \le T} \widehat{\mathcal{W}}_{j_0}^{j,(n)}(s)(-\infty, \widehat{F}_{k_0,j_0}^{(n)}(s)) = o(\sqrt{n})
\end{aligned}
$$

because of (4.7). The last term, $R^{(n)}(nt)$, satisfies

$$(4.16) \qquad R^{(n)}(nt) \leq \max_{0 \leq s \leq T} [W_{j_0}^{(n)}(ns) - W_{j_0}^{(n)}(ns-)] = o(\sqrt{n})$$

as in (4.14). We also have

$$
\begin{aligned}
K_k^{(n)}(nt) &= \sum_{i=1}^{A_{k,j_0}^{(n)}(nt)} \left( v_{k,j_0}^{(n)} - \frac{1}{\mu_{k,j_0}^{(n)}} \right) - \sum_{i=1}^{A_{k,j_0}^{(n)}(n\tau_{k_0,j_0}^{(n)}(t))} \left( v_{k,j_0}^{(n)} - \frac{1}{\mu_{k,j_0}^{(n)}} \right) \\
&\quad + \frac{1}{\mu_{k,j_0}^{(n)}} [A_{k,j_0}^{(n)}(nt) - A_{k,j_0}^{(n)}(n\tau_{k_0,j_0}^{(n)}(t))] \\
&= \sqrt{n} \left[ \widehat{V}_{k,j_0}^{(n)} \left( \frac{1}{n} A_{k,j_0}^{(n)}(nt) \right) - \widehat{V}_{k,j_0}^{(n)} \left( \frac{1}{n} A_{k,j_0}^{(n)}(n\tau_{k_0,j_0}^{(n)}(t)) \right) \right] \\
&\quad + \frac{\sqrt{n}}{\mu_{k,j_0}^{(n)}} [\widehat{A}_{k,j_0}^{(n)}(t) - \widehat{A}_{k,j_0}^{(n)}(\tau_{k_0,j_0}^{(n)}(t))] + \rho_{k,j_0}^{(n)} n(t - \tau_{k_0,j_0}^{(n)}(t)) \\
&= \sqrt{n} \Bigg[ \widehat{V}_{k,j_0}^{(n)} \left( \frac{1}{\sqrt{n}} \widehat{A}_{k,j_0}^{(n)}(t) + \lambda_k t \right) \\
&\qquad\quad - \widehat{V}_{k,j_0}^{(n)} \left( \frac{1}{\sqrt{n}} \widehat{A}_{k,j_0}^{(n)}(\tau_{k_0,j_0}^{(n)}(t)) + \lambda_k \tau_{k_0,j_0}^{(n)}(t) \right) \Bigg] \\
&\quad + \frac{\sqrt{n}}{\mu_{k,j_0}^{(n)}} [\widehat{A}_{k,j_0}^{(n)}(t) - \widehat{A}_{k,j_0}^{(n)}(\tau_{k_0,j_0}^{(n)}(t))] + \rho_{k,j_0}^{(n)} n(t - \tau_{k_0,j_0}^{(n)}(t)).
\end{aligned}
$$

(4.17)

From (2.21) and (2.24), we obtain

$$(4.18) \qquad K_k^{(n)}(nt) = \rho_{k,j_0}^{(n)} n(t - \tau_{k_0,j_0}^{(n)}(t)) + O(\sqrt{n}).$$

Finally, we estimate $D^{(n)}(nt)$. For this we choose $y \leq y_{k_0}^*$ and divide the analysis into the two cases

$$(4.19) \qquad t - \tau_{k_0,j_0}^{(n)}(t) < \frac{1}{\sqrt{n}}(y_{k_0}^* - y)$$

and the complementary case

$$(4.20) \qquad n\tau_{k_0,j_0}^{(n)}(t) + \sqrt{n}(y_{k_0}^* - y) \leq nt.$$

In the former case, $t - \tau_{k_0,j_0}^{(n)}(t) = O(n^{-1/2})$. We show this is also true in the latter case. Under condition (4.20), because

$$F_{k_0,j_0}^{(n)}(nt) = F_{k_0,j_0}^{(n)}(n\tau_{k_0,j_0}^{(n)}(t)) - n(t - \tau_{k_0,j_0}^{(n)}(t)) \leq \sqrt{n} y_{k_0}^* - n(t - \tau_{k_0,j_0}^{(n)}(t)),$$

we have

$$D^{(n)}(nt)$$

$$\leq \sum_{i=1}^{\infty} v_{0,k_0,j_0}^{i,(n)} \mathbb{1}_{\{n\tau_{k_0,j_0}^{(n)}(t) < S_{k_0}^{i,(n)} \leq n\tau_{k_0,j_0}^{(n)}(t) + \sqrt{n}(y_{k_0}^* - y)\}}$$

$$+ \sum_{i=1}^{\infty} v_{0,k_0,j_0}^{i,(n)} \mathbb{I}_{\{n\tau_{k_0,j_0}^{(n)}(t)+\sqrt{n}(y_{k_0}^*-y)<S_{k_0}^{i,(n)}\le nt\}} \mathbb{I}_{\{L_{k_0}^{i,(n)}\le\sqrt{n}y\}}$$

$$= V_{0,k_0,j_0}^{(n)}(A_{k_0}(n\tau_{k_0,j_0}^{(n)}(t)+\sqrt{n}(y_{k_0}^*-y))) - V_{0,k_0,j_0}^{(n)}(A_{k_0}(n\tau_{k_0,j_0}^{(n)}(t)))$$

$$+ \sum_{i=1}^{A_{k_0}^{(n)}(nt)}\left[v_{0,k_0,j_0}^{i,(n)}\mathbb{I}_{\{L_{k_0}^{i,(n)}\le\sqrt{n}y\}} - \frac{1}{\mu_{k_0,j_0}^{(n)}}G_{k_0}(y)\right]$$

$$- \sum_{i=1}^{A_{k_0}^{(n)}(n\tau_{k_0,j_0}^{(n)}(t)+\sqrt{n}(y_{k_0}^*-y))}\left[v_{0,k_0,j_0}^{i,(n)}\mathbb{I}_{\{L_{k_0}^{i,(n)}\le\sqrt{n}y\}} - \frac{1}{\mu_{k_0,j_0}^{(n)}}G_{k_0}(y)\right]$$

(4.21)
$$+ \frac{G_{k_0}(y)}{\mu_{k_0,j_0}^{(n)}}[A_{k_0}^{(n)}(nt) - A_{k_0}^{(n)}(n\tau_{k_0,j_0}^{(n)}(t)+\sqrt{n}(y_{k_0}^*-y))]$$

$$= \sqrt{n}\left[\widehat{M}_{0,k_0,j_0}^{(n)}\left(\tau_{k_0,j_0}^{(n)}(t)+\frac{1}{\sqrt{n}}(y_{k_0}^*-y)\right) - \widehat{M}_{0,k_0,j_0}^{(n)}(\tau_{k_0,j_0}^{(n)}(t))\right]$$

$$+ \sqrt{n}\left[\widehat{T}_{0,k_0,j_0}^{(n)}\left(\frac{1}{n}A_{k_0}^{(n)}(nt);y\right)\right.$$

$$\left. - \widehat{T}_{0,k_0,j_0}^{(n)}\left(\frac{1}{n}A_{k_0}^{(n)}(n\tau_{k_0,j_0}^{(n)}(t)+\sqrt{n}(y_{k_0}^*-y));y\right)\right]$$

$$+ \frac{\sqrt{n}G_{k_0}(y)}{\mu_{k_0,j_0}^{(n)}}\left[\widehat{A}_{k_0}^{(n)}(t) - \widehat{A}_{k_0}^{(n)}\left(\tau_{k_0,j_0}^{(n)}(t)+\frac{1}{\sqrt{n}}(y_{k_0}^*-y)\right)\right]$$

$$+ nG_{k_0}(y)\rho_{k_0,j_0}^{(n)}(t-\tau_{k_0,j_0}^{(n)}(t)) + \sqrt{n}(1-G_{k_0}(y))\rho_{k_0,j_0}^{(n)}(y_{k_0}^*-y).$$

From (2.20), (2.22) and (2.23), we obtain

(4.22)
$$D^{(n)}(nt) \le nG_{k_0}(y)\rho_{k_0,j_0}^{(n)}(t-\tau_{k_0,j_0}^{(n)}(t))$$
$$+ \sqrt{n}(1-G_{k_0}(y))\rho_{k_0,j_0}^{(n)}(y_{k_0}^*-y) + O(\sqrt{n}).$$

Substitution of (4.14)–(4.16), (4.18) and (4.22) into (4.13) yields

(4.23)
$$0 \le \mathcal{W}_{j_0}^{(n)}(nt)(-\infty, F_{k_0,j_0}^{(n)}(nt))$$
$$\le n(t-\tau_{k_0,j_0}^{(n)}(t))\left[\sum_{k\in\mathcal{C}(j_0)}\rho_{k,j_0}^{(n)} - 1 - (1-G_{k_0}(y))\rho_{k_0,j_0}^{(n)}\right] + O(\sqrt{n})$$
$$\le -n(t-\tau_{k_0,j_0}^{(n)}(t))(1-G_{k_0}(y))\rho_{k_0,j_0}^{(n)} + O(\sqrt{n}),$$

where the last inequality follows from (2.17). Assume for the moment that $y < y_{k_0}^*$. Then $(1-G_{k_0}(y))\rho_{k_0,j_0}^{(n)}$ is strictly positive and bounded away from zero uniformly in $n$. This implies

$$t - \tau_{k_0,j_0}^{(n)}(t) = O(n^{-1/2}).$$

Hence, (4.11) holds.

Armed with (4.11), we return to the weaker assumption $y \leq y_{k_0}^*$ and use the differencing theorem in (4.17) (see, e.g., Theorem A.3 of [5]), (2.21) and (2.24) to obtain

$$(4.18') \qquad K_k^{(n)}(nt) = \rho_{k,j_0}^{(n)} n(t - \tau_{k_0,j_0}^{(n)}(t)) + o(\sqrt{n})$$

in place of (4.18). Similarly, (4.22) becomes

$$(4.22') \qquad \begin{aligned} D^{(n)}(nt) &\leq n G_{k_0}(y) \rho_{k_0,j_0}^{(n)}(t - \tau_{k_0,j_0}^{(n)}(t)) \\ &\quad + \sqrt{n}(1 - G_{k_0}(y)) \rho_{k_0,j_0}^{(n)}(y_{k_0}^* - y) + o(\sqrt{n}). \end{aligned}$$

If (4.20) holds, we may substitute (4.14)–(4.16), (4.18') and (4.22') into (4.13) to obtain

$$(4.24) \qquad \begin{aligned} 0 &\leq \mathcal{W}_{j_0}^{(n)}(nt)(-\infty, F_{k_0,j_0}(nt)) \\ &\leq -n(t - \tau_{k_0,j_0}^{(n)}(t))\left(1 - G_{k_0}(y)\rho_{k_0,j_0}^{(n)} - \sum_{k \in \mathcal{C}(j_0), k \neq k_0} \rho_{k,j_0}^{(n)}\right) \\ &\quad + \sqrt{n}(1 - G_{k_0}(y))\rho_{k_0,j_0}^{(n)}(y_{k_0}^* - y) + o(\sqrt{n}). \end{aligned}$$

Assume for the moment that $y < y_{k_0}^*$. Then (4.24) implies

$$(4.25) \qquad \sqrt{n}(t - \tau_{k_0,j_0}^{(n)}(t)) \leq C_n(y_{k_0}^* - y) + o(1),$$

where

$$C_n = \frac{(1 - G_{k_0}(y))\rho_{k_0,j_0}^{(n)}}{1 - G_{k_0}(y)\rho_{k_0,j_0}^{(n)} - \sum_{k \in \mathcal{C}(j_0), k \neq k_0} \rho_{k,j_0}^{(n)}}.$$

The constants $C_n$ converge to a finite limit as $n \to \infty$, which is bounded uniformly in $y < y_{k_0}^*$. Since $y$ may be arbitrarily close to $y_{k_0}^*$, by (4.25), if (4.20) holds, then (4.12) holds also. Similarly, if (4.19) holds, then we have (4.25) with $C_n \equiv 1$ and again (4.12) follows.

From (4.12) we immediately obtain (4.5) for $k = k_0$ and $j = j_0$, because

$$\begin{aligned} &\mathcal{Q}_{j_0}^{(n)}(nt)(-\infty, F_{k_0,j_0}^{(n)}(nt)) \\ &\leq \sum_{k \in \mathcal{C}(j_0)} [A_{k,j_0}^{(n)}(nt) - A_{k,j_0}^{(n)}(n\tau_{k_0,j_0}^{(n)}(t)-)] \\ &= \sqrt{n} \sum_{k \in \mathcal{C}(j_0)} [\widehat{A}_{k,j_0}^{(n)}(t) - \widehat{A}_{k,j_0}^{(n)}(\tau_{k_0,j_0}^{(n)}(t)-)] + n(t - \tau_{k_0,j_0}^{(n)}(t)) \sum_{k \in \mathcal{C}(j_0)} \lambda_k^{(n)}. \end{aligned}$$

We divide this by $\sqrt{n}$ and use the differencing theorem, (2.24) and (4.12) to obtain the first relation in (4.5). For the second part of (4.5), we set $y = y_{k_0}^*$, so that case (4.19) is vacuous. Then (4.20) implies (4.24), which we divide by $\sqrt{n}$. The conclusion follows.

It remains to prove (4.6) for $k = k_0$, $j = j_0$ and $\ell$ satisfying $j_0 \in \mathcal{S}(k_0|\ell)$. From (4.8) and (4.10) we see that all work at station $j_0$ for station $\ell$ present and having lead time in $(-\infty, F_{k_0,j_0}^{(n)}(nt))$ at time $nt$ must arrive in the time interval $[n\tau_{k_0,j_0}^{(n)}(t), nt]$. It follows that

$$(4.26) \quad \mathcal{W}_\ell^{j_0,(n)}(nt)(-\infty, F_{k_0,j_0}^{(n)}(nt)) \le H_\ell^{(n)}(nt) + \sum_{k \in \mathcal{C}(j_0) \cap \mathcal{C}(\ell)} K_{k,\ell}^{(n)}(nt),$$

where

$$H_\ell^{(n)}(nt) \triangleq \mathcal{W}_\ell^{j_0,(n)}(n\tau_{k_0,j_0}^{(n)}(t))(-\infty, F_{k_0,j_0}^{(n)}(n\tau_{k_0,j_0}^{(n)}(t)))$$

accounts for station $\ell$ work arriving to station $j_0$ at time $n\tau_{k_0,j_0}^{(n)}(t)$, and

$$K_{k,\ell}^{(n)}(nt) \triangleq \sum_{i=1}^{\infty} v_{k,\ell}^{i,j_0,(n)} \mathbb{I}_{\{A_{k,j_0}(n\tau_{k_0,j_0}^{(n)}(t)) < i \le A_{k,j_0}^{(n)}(nt)\}}$$

is the class $k$, station $\ell$ work arriving to station $j_0$ during the time interval $(n\tau_{k_0,j_0}^{(n)}(t), nt]$. The random variables $v_{k,\ell}^{i,j_0,(n)}$, $i = 1, 2, \ldots,$ are a random permutation of the class $k$, station $\ell$ service time random variables $v_{0,k,\ell}^{i,(n)}$, $i = 1, 2, \ldots.$ The latter are indexed in order of arrival to the system; the former are indexed in order of arrival to station $j_0$. Because the index $i$ of arrival of a customer of class $k$ to station $j_0$ is independent of the service time $v_{k,\ell}^{i,j_0,(n)}$ of that customer, the sequence $v_{k,\ell}^{i,j_0,(n)}$, $i = 1, 2, \ldots,$ is independent and identically distributed, with the same distribution as $v_{0,k,\ell}^{i,(n)}$, $i = 1, 2, \ldots.$

We bound the terms appearing on the right-hand side of (4.26). We have first of all that

$$\frac{1}{\sqrt{n}} H_\ell^{(n)}(nt)$$

$$\le \sum_{k \in \mathcal{C}(\ell)} \max_{0 \le i \le A_k^{(n)}(nT)} \frac{1}{\sqrt{n}} v_{0,k,\ell}^{i,(n)}$$

$$\le \sum_{k \in \mathcal{C}(\ell)} \max_{0 \le s \le T} \left( \widehat{V}_{0,k,\ell}^{(n)}\left( \frac{1}{n} A_k^{(n)}(ns) \right) - \widehat{V}_{0,k,\ell}^{(n)}\left( \frac{1}{n} A_k^{(n)}(ns-) \right) + \frac{1}{\sqrt{n}\mu_{k,\ell}^{(n)}} \right).$$

But the process

$$\widehat{V}_{0,k,\ell}^{(n)}\left( \frac{1}{n} A_k^{(n)}(ns) \right) = \widehat{V}_{0,k,\ell}^{(n)}\left( \frac{1}{\sqrt{n}} \widehat{A}_k^{(n)}(s) + \lambda_k^{(n)} s \right)$$

converges weakly in $D_{\mathbb{R}}[0, \infty)$ to the continuous process $\widehat{V}^*_{k,\ell}(\lambda_k s)$ [see (2.21) and (2.22)] and, thus, its maximum jump over $s \in [0, T]$ converges to zero. It follows that

$$(4.27) \qquad \max_{0 \le t \le T} \frac{1}{\sqrt{n}} H^{(n)}_\ell(nt) \Rightarrow 0.$$

Let us now define

$$\widehat{V}^{j_0,(n)}_{k,\ell}(t) \triangleq \frac{1}{\sqrt{n}} \sum_{i=1}^{\lfloor nt \rfloor} \left( v^{i,j_0,(n)}_{k,\ell} - \frac{1}{\mu^{(n)}_{k,\ell}} \right),$$

which satisfies $\widehat{V}^{j_0,(n)}_{k,\ell} \Rightarrow \widehat{V}^*_{k,\ell}$, where $\widehat{V}^*_{k,\ell}$ is a continuous process [cf. (2.21)]. Then

$$\begin{aligned}
\frac{1}{\sqrt{n}} K^{(n)}_{k,\ell}(nt) &= \widehat{V}^{j_0,(n)}_{k,\ell}\left( \frac{1}{n} A^{(n)}_{k,j_0}(nt) \right) - \widehat{V}^{j_0,(n)}_{k,\ell}\left( \frac{1}{n} A^{(n)}_{k,j_0}(n\tau^{(n)}_{k_0,j_0}(t)) \right) \\
&\quad + \frac{1}{\sqrt{n}\mu^{(n)}_{k,\ell}}(A^{(n)}_{k,j_0}(nt) - A^{(n)}_{k,j_0}(n\tau^{(n)}_{k_0,j_0}(t))) \\
&= \widehat{V}^{j_0,(n)}_{k,\ell}\left( \frac{1}{\sqrt{n}} \widehat{A}^{(n)}_{k,j_0}(t) + \lambda^{(n)}_k t \right) \\
&\quad - \widehat{V}^{j_0,(n)}_{k,\ell}\left( \frac{1}{\sqrt{n}} \widehat{A}^{(n)}_{k,j_0}(\tau^{(n)}_{k_0,j_0}(t)) + \lambda^{(n)}_k \tau^{(n)}_{k_0,j_0}(t) \right) \\
&\quad + \frac{1}{\mu^{(n)}_{k,\ell}}(\widehat{A}^{(n)}_{k,j_0}(t) - \widehat{A}^{(n)}_{k,j_0}(\tau^{(n)}_{k_0,j_0}(t))) + \rho^{(n)}_{k,j_0}\sqrt{n}(t - \tau^{(n)}_{k_0,j_0}(t)).
\end{aligned}$$

The right-hand side converges weakly to zero in $D_{\mathbb{R}}[0, \infty)$ because of (4.11), (4.12), the continuity of $\widehat{V}^*_{k,\ell}$ and Assumption 2.1. In particular,

$$(4.28) \qquad \max_{0 \le t \le T} \frac{1}{\sqrt{n}} K^{(n)}_{k,\ell}(nt) \Rightarrow 0 \qquad \forall k \in \mathcal{C}(j_0) \cap \mathcal{C}(\ell).$$

From (4.26)–(4.28), we have

$$\max_{0 \le t \le T} \widehat{\mathcal{W}}^{j_0,(n)}_\ell(t)(-\infty, \widehat{F}^{(n)}_{k_0,j_0}(t)) = \max_{0 \le t \le T} \frac{1}{\sqrt{n}} \mathcal{W}^{j_0,(n)}_\ell(nt)(-\infty, F^{(n)}_{k_0,j_0}(nt)) \Rightarrow 0.$$

This gives us (4.6) for $k = k_0$, $j = j_0$ and $\ell$ satisfying $j_0 \in \mathcal{S}(k_0 | \ell)$. $\quad\square$

COROLLARY 4.4 (Crushing). *For every $j$, we have*

$$(4.29) \qquad \widehat{\mathcal{Q}}^{(n)}_j(-\infty, \widehat{F}^{(n)}_j(t)) \Rightarrow 0, \qquad \widehat{\mathcal{W}}^{(n)}_j(-\infty, \widehat{F}^{(n)}_j(t)) \Rightarrow 0.$$

Proof. From (4.5) we have

$$\widehat{\mathcal{Q}}_j^{(n)}(-\infty, \widehat{F}_j^{(n)}(t)) = \max_{k \in \mathcal{C}(j)} \widehat{\mathcal{Q}}_j^{(n)}(-\infty, \widehat{F}_{k,j}^{(n)}(t)) \Rightarrow 0,$$

$$\widehat{\mathcal{W}}_j^{(n)}(-\infty, \widehat{F}_j^{(n)}(t)) = \max_{k \in \mathcal{C}(j)} \widehat{\mathcal{W}}_j^{(n)}(-\infty, \widehat{F}_{k,j}^{(n)}(t)) \Rightarrow 0. \qquad \square$$

COROLLARY 4.5. *For all stations $j$ and $\ell$ for which there exists a customer class $k$ such that $j \in \mathcal{S}(k|\ell)$, we have*

$$\widehat{\mathcal{W}}_\ell^{j,(n)}(-\infty, \widehat{F}_j^{(n)}) \Rightarrow 0.$$

Proof. It suffices to show that for every $k \in \mathcal{C}(j) \cap \mathcal{C}(\ell)$, we have

$$\widehat{\mathcal{W}}_{k,\ell}^{j,(n)}(-\infty, \widehat{F}_j^{(n)}) \Rightarrow 0.$$

For such a $k$, let $\mathcal{I}_{k,j}^{(n)}(nt)$ denote the indices $i$, according to the order of arrival to the system, of the class $k$ customers at station $j$ at time $nt$ with lead times in $(-\infty, F_j^{(n)}(nt))$. Then

$$\widehat{\mathcal{W}}_{k,\ell}^{j,(n)}(t)(-\infty, \widehat{F}_j^{(n)}(t)) \leq \frac{1}{\sqrt{n}} \sum_{i \in \mathcal{I}_{k,j}^{(n)}(nt)} v_{0,k,\ell}^{i,(n)}.$$

Let $|\mathcal{I}_{k,j}^{(n)}(nt)|$ denote the cardinality of $\mathcal{I}_{k,j}^{(n)}(nt)$. For each positive integer $m$, let $\alpha_m^{(n)} = \mathbb{P}\{|\mathcal{I}_{k,j}^{(n)}(nt)| = m\}$. Let $\mathbb{P}_0^{(n)}$ be the zero measure, and for each positive integer $m$, let $\mathbb{P}_m^{(n)}$ denote the measure induced on $\mathbb{R}$ by the random variable $\frac{1}{m} \sum_{i \in \mathcal{I}_{k,j}^{(n)}(nt)} v_{0,k,\ell}^{i,(n)}$. Conditioned on $|\mathcal{I}_{k,j}^{(n)}(nt)| = m$, the distribution of $\sum_{i \in \mathcal{I}_{k,j}^{(n)}(nt)} v_{0,k,\ell}^{i,(n)}$ is the same as the distribution of $\sum_{i=1}^m v_{0,k,\ell}^{i,(n)}$. Because the random variables $v_{0,k,\ell}^{i,(n)}$ do not enter the determination of the indices which below to $\mathcal{I}_{k,j}^{(n)}$, the measure induced on $\mathbb{R}$ by $\frac{1}{|\mathcal{I}_{k,j}^{(n)}|} \mathbb{I}_{\{|\mathcal{I}_{k,j}^{(n)}(nt)| \geq 1\}} \sum_{i \in \mathcal{I}_{k,j}^{(n)}(nt)} v_{0,k,\ell}^{i,(n)}$ is $\sum_{n=0}^\infty \alpha_m^{(n)} \mathbb{P}_m^{(n)}$. The set of probability measure $\{\mathbb{P}_m^{(n)}\}$ is tight because for $m \geq 1$ and $K \geq (\mu_{k,\ell}^{(n)})^{-1}$, we have

$$\mathbb{P}\left\{\frac{1}{m}\sum_{i=1}^m v_{0,k,\ell}^{i,(n)} \geq K\right\} \leq \mathbb{P}\left\{\frac{1}{m}\sum_{i=1}^m (v_{0,k,\ell}^{i,(n)} - (\mu_{k,\ell}^{(n)})^{-1}) \geq K - (\mu_{k,\ell}^{(n)})^{-1}\right\}$$

$$\leq \frac{1}{(K - (\mu_{k,\ell}^{(n)})^{-1})^2} \mathbb{E}\left[\frac{1}{m}\sum_{i=1}^m (v_{0,k,\ell}^{i,(n)} - (\mu_{k,\ell}^{(n)})^{-1})\right]^2$$

$$= \frac{(\beta_{k,\ell}^{(n)})^2}{m(K - (\mu_{k,\ell}^{(n)})^{-1})^2},$$

and this can be made arbitrarily small, uniformly in $m$ and $n$, by the choice of $K$. Consequently, the set of probability measures $\{\sum_{m=0}^{\infty} \alpha_m^{(n)} \mathbb{P}_m^{(n)}\}_{n \geq 1}$ is also tight. In particular, given $\varepsilon > 0$, there exists $K > 0$ such that

$$\mathbb{P}\left\{\frac{1}{|\mathcal{J}_{k,j}^{(n)}(nt)|} \mathbb{I}_{\{|\mathcal{J}_{k,j}^{(n)}(nt)| \geq 1\}} \sum_{i \in \mathcal{J}_{k,j}^{(n)}(nt)} v_{0,k,\ell}^{(n)} \leq K\right\} \geq 1 - \varepsilon$$

for all $n \geq 1$.

According to the first part of (4.29), $\frac{1}{\sqrt{n}}|\mathcal{J}_{k,j}^{(n)}(nt)| \Rightarrow 0$, and hence there is an integer $N$ such that

$$\mathbb{P}\left\{\frac{1}{\sqrt{n}}|\mathcal{J}_{k,j}^{(n)}(nt)| \leq \frac{\varepsilon}{K}\right\} \geq 1 - \varepsilon$$

for all $n \geq N$. Therefore, for $n \geq N$,

$$\mathbb{P}\{\widehat{W}_{k,\ell}^{j,(n)}(t)(-\infty, \widehat{F}_j^{(n)}(t)) \leq \varepsilon\}$$

$$\geq \mathbb{P}\left\{\frac{1}{|\mathcal{J}_{k,j}^{(n)}(nt)|} \mathbb{I}_{\{|\mathcal{J}_{k,j}^{(n)}(nt)| \geq 1\}} \sum_{i \in \mathcal{J}_{k,j}^{(n)}(nt)} v_{0,k,j}^{i,(n)} \cdot \frac{1}{\sqrt{n}}|\mathcal{J}_{k,j}^{(n)}(nt)| \leq \varepsilon\right\}$$

$$\geq 1 - 2\varepsilon.$$

This establishes the corollary.    □

The following lemma gives a tightness bound for the scaled frontiers.

LEMMA 4.6.   *For every $T > 0$, $\varepsilon > 0$, $j \in \{1, \dots, J\}$ and $k \in \mathcal{C}(j)$, there exists $y \in (-\infty, y_k^*)$ such that for all $n$,*

$$(4.30) \qquad \mathbb{P}\left\{\inf_{0 \leq t \leq T} \widehat{F}_{k,j}^{(n)}(t) < y\right\} < \varepsilon.$$

PROOF.   As in the proof of Lemma 4.3, we fix $k$ and proceed by induction along $\mathcal{P}(k)$. Let $\ell \in \mathcal{P}(k)$ be given and assume that for $j \in \mathcal{S}(k|\ell), T > 0$ and $\varepsilon > 0$, the corresponding $y$ satisfying (4.30) for all $n$ can be found. In particular, no assumption is necessary to analyze the first station in $\mathcal{P}(k)$.

We first argue that

$$(4.31) \qquad \widehat{W}_\ell^{(n)}(t) \geq \widehat{\mathcal{V}}_{0,k,\ell}^{(n)}(t)\left(\widehat{F}_{k,\ell}^{(n)}(t), \min_{j \in \mathcal{S}(k|\ell)} \widehat{F}_{k,j}^{(n)}(t)\right) + o(1).$$

Indeed, the workload at station $\ell$ is at least as great as the workload brought to station $\ell$ by class $k$ customers with lead times in $(F_{k,\ell}^{(n)}(nt), \min_{j \in \mathcal{S}(k|\ell)} F_{k,j}^{(n)}(nt))$. None of class $k$ customers who have arrived to the system by time $nt$ and have lead times at this time greater than $F_{k,\ell}^{(n)}(nt)$ has ever been in service at station $\ell$ by time $nt$. Thus, every such customer is either in queue at station $\ell$ or in queue at some station $j_0 \in \mathcal{S}(k|\ell)$. By Lemma 4.3, for such $j_0$ we have

$$\widehat{\mathcal{W}}_\ell^{j_0,(n)}(t)\left(-\infty, \min_{j \in \mathcal{S}(k|\ell)} \widehat{F}_{k,j}^{(n)}(t)\right) \leq \widehat{\mathcal{W}}_\ell^{j_0,(n)}(-\infty, \widehat{F}_{k,j_0}^{(n)}(t)) = o(1),$$

so the difference between $\widehat{\mathcal{V}}_{0,k,\ell}^{(n)}(\widehat{F}_{k,\ell}^{(n)}(t), \min_{j \in \mathcal{S}(k|\ell)} \widehat{F}_{k,j}^{(n)}(t))$ and the scaled workload for station $\ell$ associated with class $k$ customers already present at $\ell$ with lead times in $(F_{k,\ell}^{(n)}(nt), \min_{j \in \mathcal{S}(k|\ell)} F_{k,j}^{(n)}(nt))$ is of the order $o(1)$. This justifies (4.31).

Fix $T > 0$ and $\varepsilon > 0$. By the induction hypothesis, there exists $y_1 < y_k^*$ such that, for all $n$, $\mathbb{P}(A_n) \geq 1 - \frac{\varepsilon}{4}$, where

$$A_n \triangleq \left\{ \inf_{0 \leq t \leq T} \min_{j \in \mathcal{S}(k|\ell)} \widehat{F}_{k,j}^{(n)}(t) > y_1 \right\}.$$

By (2.25) and the continuous mapping theorem,

$$\sup_{0 \leq t \leq T} \widehat{W}_\ell^{(n)}(t) \Rightarrow \sup_{0 \leq t \leq T} W_\ell^*(t).$$

Inequality (4.31) and the fact that $\lim_{y \to -\infty} H_k(y) = \infty$ enables us to choose $y_2 < y_1$ such that, for all $n$, $\mathbb{P}(B_n) \geq 1 - \frac{\varepsilon}{4}$, where

$$B_n \triangleq \left\{ \sup_{0 \leq t \leq T} \widehat{\mathcal{V}}_{0,k,\ell}^{(n)}(t)\left(\widehat{F}_{k,\ell}^{(n)}(t), \min_{j \in \mathcal{S}(k|\ell)} \widehat{F}_{k,j}^{(n)}(t)\right) \leq \sqrt{H_k(y_2)} \right\}.$$

On $A_n \cap B_n$, for $0 \leq t \leq T$, we have

$$(4.32) \qquad \sup_{0 \leq t \leq T} \widehat{\mathcal{V}}_{0,k,\ell}^{(n)}(t)(\widehat{F}_{k,\ell}^{(n)}(t), y_1) \leq \sqrt{H_k(y_2)}.$$

By Proposition 4.1, we can find $N$ such that for all $n \geq N$, $\mathbb{P}(C_n) \geq 1 - \frac{\varepsilon}{4}$, where

$$C_n \triangleq \left\{ \sup_{y_2 \leq y \leq y_k^*} \sup_{0 \leq t \leq T} \left| \widehat{\mathcal{V}}_{0,k,\ell}^{(n)}(t)(y, \infty) + \rho_{k,\ell}[H_k(y + \sqrt{n}t) - H_k(y)] \right| \right.$$
$$\left. \leq \frac{\rho_{k,\ell}}{2} H_k(y_1) \right\}.$$

For $n \geq N$, $\mathbb{P}(A_n \cap B_n \cap C_n) \geq 1 - \frac{3\varepsilon}{4}$. By (4.32), on $A_n \cap B_n \cap C_n$ we have

$$\sup_{0 \leq t \leq T} \widehat{\mathcal{V}}_{0,k,\ell}^{(n)}(t)(y_2, y_1] \mathbb{I}_{\{\widehat{F}_{k,\ell}^{(n)}(t) < y_2\}} \leq \sqrt{H_k(y_2)},$$

so

$$
\begin{aligned}
\sqrt{H_k(y_2)} &+ \sup_{0 \le t \le T} \widehat{\mathcal{V}}_{0,k,\ell}^{(n)}(t)(y_1,\infty) \mathbb{I}_{\{\widehat{F}_{k,\ell}^{(n)}(t) < y_2\}} \\
&\ge \sup_{0 \le t \le T} \widehat{\mathcal{V}}_{0,k,\ell}^{(n)}(t)(y_2,\infty) \mathbb{I}_{\{\widehat{F}_{k,\ell}^{(n)}(t) < y_2\}} \\
&\ge \frac{\rho_{k,\ell}}{2} H_k(y_2) \max_{0 \le t \le T} \mathbb{I}_{\{\widehat{F}_{k,\ell}^{(n)}(t) < y_2\}} \\
&= \frac{\rho_{k,\ell}}{2} H_k(y_2) \mathbb{I}_{\{\inf_{0 \le t \le T} \widehat{F}_{k,\ell}^{(n)}(t) < y_2\}},
\end{aligned}
$$

$(4.33)$

where the third line follows from the definition of $C_n$, the fact that $y_2 < y_1 < y_k^*$ implies $H_k(y_2) > H_k(y_1)$, and the inequality $\widehat{F}_{k,\ell}^{(n)}(t) \ge y_k^* - \sqrt{n}\,t$ (following immediately from the definition of the frontier) resulting in

$$
(4.34) \qquad y_k^* < y_2 + \sqrt{n}\,t \text{ on } \{\widehat{F}_{k,\ell}^{(n)}(t) < y_2\}.
$$

Also, on $A_n \cap B_n \cap C_n$,

$$
\begin{aligned}
\sup_{0 \le t \le T} &\widehat{\mathcal{V}}_{0,k,\ell}^{(n)}(t)(y_1,\infty) \mathbb{I}_{\{\widehat{F}_{k,\ell}^{(n)}(t) < y_2\}} \\
&= \sup_{0 \le t \le T} [\widehat{\mathcal{V}}_{0,k,\ell}^{(n)}(t)(y_1,\infty) + \rho_{k,\ell} H_k(y_1 + \sqrt{n}\,t)] \mathbb{I}_{\{\widehat{F}_{k,\ell}^{(n)}(t) < y_2\}} \\
&\le \frac{3\rho_{k,\ell}}{2} H_k(y_1) \triangleq c.
\end{aligned}
$$

[The second line follows from $y_2 < y_1$ and (4.34) and the third one from the definition of $C_n$.] Thus, (4.33) yields, for $n \ge N$, on $A_n \cap B_n \cap C_n$,

$$
\frac{\rho_{k,\ell}}{2} H_k(y_2) \mathbb{I}_{\{\inf_{0 \le t \le T} \widehat{F}_{k,\ell}^{(n)}(t) < y_2\}} \le \sqrt{H_k(y_2)} + c,
$$

and, therefore,

$$
\begin{aligned}
\mathbb{P}\Big\{ &\inf_{0 \le t \le T} \widehat{F}_{k,\ell}^{(n)}(t) < y_2 \Big\} - \frac{3\varepsilon}{4} \\
&\le \mathbb{P}\Big( \Big\{ \inf_{0 \le t \le T} \widehat{F}_{k,\ell}^{(n)}(t) < y_2 \Big\} \cap A_n \cap B_n \cap C_n \Big) \\
&\le \frac{2(\sqrt{H_k(y_2)} + c)}{\rho_{k,\ell} H_k(y_2)} < \frac{\varepsilon}{4}
\end{aligned}
$$

for $y_2$ small enough. Thus, (4.30) holds for $k$, $j = \ell$, $y = y_2$ and $n \ge N$. Taking $y$ smaller, if necessary, we extend (4.30) to $k$, $j = \ell$ and all $n$.  $\square$

COROLLARY 4.7.  *For every $j$, we have*

$$
\widehat{\mathcal{Q}}_j^{(n)}(-\infty, \widehat{F}_j^{(n)}(t)] \Rightarrow 0, \qquad \widehat{\mathcal{W}}_j^{(n)}(-\infty, \widehat{F}_j^{(n)}(t)] \Rightarrow 0.
$$

*Moreover, for all stations $j$ and $\ell$ for which there exists a customer class $k$ such that $j \in \mathcal{S}(k|\ell)$, we have*

$$\widehat{\mathcal{W}}_\ell^{j,(n)}(-\infty, \widehat{F}_j^{(n)}(t)] \Rightarrow 0.$$

PROOF. This is just a restatement of Corollaries 4.4 and 4.5, except that the half-line $(-\infty, \widehat{F}_j^{(n)}(t)]$ is now closed on the right. Corollary 4.2 asserts that if $\widehat{F}_j^{(n)}$ were bounded below, uniformly in $t \in [0, T]$ and $n$, then the inclusion of this endpoint would make no difference. Using Lemma 4.6, we can ensure that $\widehat{F}_j^{(n)}$ is bounded below with probability arbitrarily close to 1, and the result follows. □

**5. Inverting the frontier equations.** In this section we show the first part of Theorem 3.2, that is, that the function $\Phi$ defined by (3.20) is a homeomorphism of the set $D$ given by (3.12) and (3.13) onto $[0, \infty)^J$ (Proposition 5.5). It is clear that $\Phi$ is continuous. Lemma 5.1 asserts that $\Phi$ maps $D$ onto $[0, \infty)^J$. The proof of this lemma contains an explicit algorithm for inverting $\Phi$. Lemmas 5.2 and 5.3 show that $\Phi$ is one-to-one on $D$. Finally, Lemma 5.4, examining the limiting behavior of $\Phi(y)$ in $D$ as $y \to \infty$, is used to show that the mapping $\Phi$ is open.

LEMMA 5.1. $\Phi(D) = [0, \infty)^J$.

PROOF. Let $w = (w_1, \ldots, w_J) \in [0, \infty)^J$ be given. The aim is to find $y = (y_1, \ldots, y_J) \in D$ such that $\Phi(y) = w$, that is, to solve the frontier equations

$$(5.1) \quad w_j = \sum_{k \in \mathcal{C}(j)} \rho_{k,j} \left[ H_k(y_j) - H_k \left( \min_{i \in \mathcal{S}(k|j)} y_i \right) \right]^+, \qquad j = 1, \ldots, J,$$

for a $y \in D$. Note that if $k$ is a customer class entering the system at station $j$, then $\min_{i \in \mathcal{S}(k|j)} y_i = \min_{i \in \varnothing} y_i = \infty$, so

$$H_k(y_j) - H_k \left( \min_{i \in \mathcal{S}(k|j)} y_i \right) = H_k(y_j).$$

We may rewrite (5.1) in the form

$$w_j = \sum_{k \in \mathcal{K}_0(j)} \rho_{k,j} H_k(y_j) + \sum_{k \in \mathcal{C}(j) \setminus \mathcal{K}_0(j)} \rho_{k,j} \left[ H_k(y_j) - H_k \left( \min_{i \in \mathcal{S}(k|j)} y_i \right) \right]^+.$$

For $j \in \mathcal{J}_0$ and $y \in \mathbb{R}$, we define

$$K_{0,j}(y) = \sum_{k \in \mathcal{K}_0(j)} \rho_{k,j} H_k(y) = \sum_{k \in \mathcal{K}_0(j)} \rho_{k,j} \left[ H_k(y) - H_k \left( \min_{i \in \mathcal{S}(k|j)} y_i \right) \right]^+.$$

Although defined on all of $\mathbb{R}$, we shall be interested in $K_{0,j}$ restricted to a smaller set. In particular,

$$K_{0,j} : \left( -\infty, \max_{k \in \mathcal{K}_0(j)} y_k^* \right] \xrightarrow{\text{onto}} [0, \infty)$$

is strictly decreasing and has a strictly decreasing inverse

$$K_{0,j}^{-1} : [0, \infty) \xrightarrow{\text{onto}} \left( -\infty, \max_{k \in \mathcal{K}_0(j)} y_k^* \right].$$

We choose $j_1 \in \mathcal{J}_0$ so that $K_{0,j_1}^{-1}(w_{j_1}) = \max_{j \in \mathcal{J}_0} K_{0,j}^{-1}(w_j)$ and we set $y_{j_1} = K_{0,j_1}^{-1}(w_{j_1})$. Then

$$
\begin{align}
(5.2) \quad w_j &\geq K_{0,j}(y_{j_1}) \\
&= \sum_{k \in \mathcal{K}_0(j)} \rho_{k,j} \left[ H_k(y_{j_1}) - H_k \left( \min_{i \in \mathcal{S}(k|j)} y_i \right) \right]^+ \qquad \forall j \in \mathcal{J}_0,
\end{align}
$$

$$(5.3) \quad w_{j_1} = K_{0,j_1}(y_{j_1}) = \sum_{k \in \mathcal{K}_0(j_1)} \rho_{k,j_1} \left[ H_k(y_{j_1}) - H_k \left( \min_{i \in \mathcal{S}(k|j_1)} y_i \right) \right]^+.$$

*Induction hypothesis.* Suppose that for $m = 1, \ldots, M$, we have chosen distinct indices $j_1, j_2, \ldots, j_M$, have defined numbers $y_{j_1} \geq y_{j_2} \geq \cdots \geq y_{j_M}$, and have defined functions

$$K_{m-1,j}(y) = \sum_{k \in \mathcal{K}_{m-1}^{(j_1,\ldots,j_{m-1})}(j)} \rho_{k,j} \left[ H_k(y) - H_k \left( \min_{i \in \mathcal{S}(k|j)} y_i \right) \right]^+, \qquad y \in \mathbb{R},$$

for $j \in \mathcal{J}_{m-1}^{(j_1,\ldots,j_{m-1})}$. Although defined on all of $\mathbb{R}$, each function

$$K_{m-1,j} : \left( -\infty, \max_{k \in \mathcal{K}_{m-1}^{(j_1,\ldots,j_{m-1})}(j)} \left( y_k^* \wedge \min_{i \in \mathcal{S}(k|j)} y_i \right) \right] \xrightarrow{\text{onto}} [0, \infty)$$

is strictly decreasing when restricted to the indicated set and, therefore, has a strictly decreasing inverse

$$(5.4) \quad K_{m-1,j}^{-1} : [0, \infty) \xrightarrow{\text{onto}} \left( -\infty, \max_{k \in \mathcal{K}_{m-1}^{(j_1,\ldots,j_{m-1})}(j)} \left( y_k^* \wedge \min_{i \in \mathcal{S}(k|j)} y_i \right) \right].$$

Suppose further that for $m = 1, \ldots, M$, we have $j_m \in \mathcal{J}_{m-1}^{(j_1,\ldots,j_{m-1})}$ and

$$K_{m-1,j_m}^{-1}(w_{j_m}) = \max_{j \in \mathcal{J}_{m-1}^{(j_1,\ldots,j_{m-1})}} K_{m-1,j}^{-1}(w_j), \qquad y_{j_m} = K_{m-1,j_m}^{-1}(w_{j_m});$$

hence,

$$w_j \geq K_{m-1,j}(y_{j_m})$$

(5.5)
$$= \sum_{k \in \mathcal{K}_{m-1}^{(j_1,\ldots,j_{m-1})}(j)} \rho_{k,j} \left[ H_k(y_{j_m}) - H_k\left( \min_{i \in \mathcal{S}(k|j)} y_i \right) \right]^+$$

$$\forall j \in \mathcal{J}_{m-1}^{(j_1,\ldots,j_{m-1})},$$

$$w_{j_m} = K_{m-1,j_m}(y_{j_m})$$

(5.6)
$$= \sum_{k \in \mathcal{K}_{m-1}^{(j_1,\ldots,j_{m-1})}(j_m)} \rho_{k,j_m} \left[ H_k(y_{j_m}) - H_k\left( \min_{i \in \mathcal{S}(k|j_m)} y_i \right) \right]^+.$$

*Induction step.* If $M = J$, we terminate the construction. If $M < J$, we proceed to step $M + 1$ as follows. Recall that the set $\mathcal{J}_M^{(j_1,\ldots,j_M)}$ is the set of all stations $j$ not among $j_1,\ldots,j_M$ with the property that at least one customer class visits $j$ and the previous stations visited by this customer class are among the stations $j_1,\ldots,j_M$. If there were no such station $j$, then all external arrivals would be to the set of stations $\{j_1,\ldots,j_M\}$ and all customers exiting a station from this set would either exit the system or else proceed to another station in this set. In this situation, stations outside set $\{j_1,\ldots,j_M\}$ would not be connected to these stations, a situation we have ruled out by assumption. Hence, $\mathcal{J}_M^{(j_1,\ldots,j_M)}$ is nonempty.

For $j \in \mathcal{J}_M^{(j_1,\ldots,j_M)}$ and $y \in \mathbb{R}$, we define

$$K_{M,j}(y) = \sum_{k \in \mathcal{K}_M^{(j_1,\ldots,j_M)}(j)} \rho_{k,j} \left[ H_k(y) - H_k\left( \min_{i \in \mathcal{S}(k|j)} y_i \right) \right]^+.$$

Although defined on all of $\mathbb{R}$,

$$K_{M,j} : \left( -\infty, \max_{k \in \mathcal{K}_M^{(j_1,\ldots,j_M)}(j)} \left( y_k^* \wedge \min_{i \in \mathcal{S}(k|j)} y_i \right) \right] \xrightarrow{\text{onto}} [0, \infty)$$

is strictly decreasing when restricted to the indicated set and, therefore, has a strictly decreasing inverse

$$K_{M,j}^{-1} : [0, \infty) \xrightarrow{\text{onto}} \left( -\infty, \max_{k \in \mathcal{K}_M^{(j_1,\ldots,j_M)}(j)} \left( y_k^* \wedge \min_{i \in \mathcal{S}(k|j)} y_i \right) \right].$$

We choose $j_{M+1} \in \mathcal{J}_M^{(j_1,\ldots,j_M)}$ so that

$$K_{M,j_{M+1}}^{-1}(w_{j_{M+1}}) = \max_{j \in \mathcal{J}_M^{(j_1,\ldots,j_M)}} K_{M,j}^{-1}(w_j)$$

and set $y_{j_{M+1}} = K_{M,j_{M+1}}^{-1}(w_{j_{M+1}})$. Then

$$
\begin{aligned}
w_j &\geq K_{M,j}(y_{j_{M+1}}) \\
(5.7) \qquad &= \sum_{k \in \mathcal{K}_M^{(j_1,\ldots,j_M)}(j)} \rho_{k,j}\left[H_k(y_{j_{M+1}}) - H_k\left(\min_{i \in \mathcal{S}(k|j)} y_i\right)\right]^+ \\
&\hspace{6cm} \forall j \in \mathcal{J}_M^{(j_1,\ldots,j_M)},
\end{aligned}
$$

$$
\begin{aligned}
w_{j_{M+1}} &= K_{M,j_{M+1}}(y_{j_{M+1}}) \\
(5.8) \qquad &= \sum_{k \in \mathcal{K}_M^{(j_1,\ldots,j_M)}(y_{j_{M+1}})} \rho_{k,j_{M+1}}\left[H_k(y_{j_{M+1}}) - H_k\left(\min_{i \in \mathcal{S}(k|j_{M+1})} y_i\right)\right]^+.
\end{aligned}
$$

To complete the induction step it remains only to show that $y_{j_M} \geq y_{j_{M+1}}$.

We divide the analysis into two cases.

CASE I.   $j_{M+1} \in \mathcal{J}_{M-1}^{(j_1,\ldots,j_{M-1})}$.

In this case (5.5) implies that

$$
(5.9) \quad w_{j_{M+1}} \geq \sum_{k \in \mathcal{K}_{M-1}^{(j_1,\ldots,j_{M-1})}(j_{M+1})} \rho_{k,j_{M+1}}\left[H_k(y_{j_M}) - H_k\left(\min_{i \in \mathcal{S}(k|j_{M+1})} y_i\right)\right]^+.
$$

For $k \in \mathcal{K}_M^{(j_1,\ldots,j_M)}(j_{M+1}) \setminus \mathcal{K}_{M-1}^{(j_1,\ldots,j_{M-1})}(j_{M+1})$, we have $j_M \in \mathcal{S}(k|j_{M+1})$, so

$$
(5.10) \qquad \min_{i \in \mathcal{S}(k|j_{M+1})} y_i = y_{j_M}.
$$

It follows that

$$
\begin{aligned}
(5.11) \quad \sum_{k \in \mathcal{K}_M^{(j_1,\ldots,j_M)}(j_{M+1}) \setminus \mathcal{K}_{M-1}^{(j_1,\ldots,j_{M-1})}(j_{M+1})} \rho_{k,j_{M+1}}&\left[H_k(y_{j_M})\right. \\
&\left. - H_k\left(\min_{i \in \mathcal{S}(k|j_{M+1})} y_i\right)\right]^+ = 0.
\end{aligned}
$$

Summing (5.9) and (5.11), we obtain

$$
w_{j_{M+1}} \geq K_{M,j_{M+1}}(y_{j_M}),
$$

and, hence,

$$
y_{j_M} \geq K_{M,j_{M+1}}^{-1}(w_{j_{M+1}}) = y_{j_{M+1}}.
$$

CASE II.   $j_{M+1} \in \mathcal{J}_M^{(j_1,\ldots,j_M)} \setminus \mathcal{J}_{M-1}^{(j_1,\ldots,j_{M-1})}$.

In this case $\mathcal{K}_M^{(j_1,\ldots,j_M)}(j_{M+1}) \neq \varnothing$ but $\mathcal{K}_{M-1}^{(j_1,\ldots,j_{M-1})}(j_{M+1}) = \varnothing$. This implies that for every $k \in \mathcal{K}_M^{(j_1,\ldots,j_M)}(j_{M+1})$, we must have $j_M \in \mathcal{S}(k|j_{M+1})$. Hence, for every $k \in \mathcal{K}_M^{(j_1,\ldots,j_M)}(j_{M+1})$, equation (5.10) holds. Equation (5.8) becomes

$$w_{j_{M+1}} = \sum_{k \in \mathcal{K}_M^{(j_1,\ldots,j_M)}(j_{M+1})} \rho_{k,j_{M+1}}[H_k(y_{j_{M+1}}) - H_k(y_{j_M})]^+.$$

If $w_{j_{M+1}} > 0$, then $y_{j_{M+1}} < y_{j_M}$. If $w_{j_{M+1}} = 0$, we have by definition

$$y_{j_{M+1}} = K_{M,j_{M+1}}^{-1}(0)$$

$$= \max_{k \in \mathcal{K}_M^{(j_1,\ldots,j_M)}(j_{M+1})} \left(y_k^* \wedge \min_{i \in \mathcal{S}(k|j_{M+1})} y_i\right)$$

$$= \max_{k \in \mathcal{K}_M^{(j_1,\ldots,j_M)}(j_{M+1})} (y_k^* \wedge y_{j_M}) \leq y_{j_M},$$

where the third equality follows from (5.10). The induction step is complete.

When this construction terminates with $M = J$, we have chosen $j_1,\ldots,j_J$, a permutation of $1,\ldots,J$, and we have defined numbers $y_{j_1} \geq y_{j_2} \geq \cdots \geq y_{j_J}$, such that (5.6) holds for $m = 1,\ldots,J$. Let $j \in \{1,\ldots,J\}$ be given, and choose $m$ so that $j = j_m$. For $k \in \mathcal{C}(j) \setminus \mathcal{K}_{m-1}^{(j_1,\ldots,j_{m-1})}(j)$, the set $\mathcal{S}(k|j)$ is not a subset of $\{j_1,\ldots,j_{m-1}\}$ and, hence, $\min_{i \in \mathcal{S}(k|j)} y_i \leq y_{j_m} = y_j$. It follows that

$$(5.12) \qquad \sum_{k \in \mathcal{C}(j) \setminus \mathcal{K}_{m-1}^{(j_1,\ldots,j_{m-1})}(j)} \rho_{k,j}\left[H_k(y_j) - H_k\left(\min_{i \in \mathcal{S}(k|j)} y_i\right)\right]^+ = 0.$$

Summing (5.6) with $j_m = j$ and (5.12), we obtain (5.1).

It remains to show that $y = (y_1,\ldots,y_J) \in D$. Let $\pi = (j_1,\ldots,j_J)$. By construction, for $m = 1,\ldots,J$, $j_m \in \mathcal{J}_{m-1}^{(j_1,\ldots,j_{m-1})} = \mathcal{J}_{m-1}^\pi$, so $\pi \in \Pi$. Moreover, $y_{j_1} \geq y_{j_2} \geq \cdots \geq y_{j_J}$ and, by (5.4),

$$y_{j_m} = K_{m-1,j_m}^{-1}(w_{j_m})$$

$$\leq \max_{k \in \mathcal{K}_{m-1}^{(j_1,\ldots,j_{m-1})}(j_m)} \left(y_k^* \wedge \min_{i \in \mathcal{S}(k|j_m)} y_i\right)$$

$$\leq \max_{k \in \mathcal{K}_{m-1}^\pi(j_m)} y_k^*.$$

Thus, $y \in D^\pi$ and, hence, $y \in D$. $\quad\square$

LEMMA 5.2. *Let $w = (w_1, \ldots, w_J) \in [0, \infty)^J$ be given. Let $y = (y_1, \ldots, y_J)$ be the solution to (5.1) constructed in Lemma 5.1 and let $\tilde{y} = (\tilde{y}_1, \ldots, \tilde{y}_J)$ be another solution to (5.1). Then, for $i = 1, \ldots, J$, we have*

$$(5.13) \qquad\qquad y_i \leq \tilde{y}_i.$$

PROOF. The proof proceeds by induction. Namely, let $\pi = (j_1, \ldots, j_J)$ be the permutation constructed in the proof of Lemma 5.1. We assume that (5.13) holds for $i = j_1, \ldots, j_M$ with some $M < J$ (in particular, for $M = 0$, no assumption is needed). We want to show that (5.13) holds for $i = j_{M+1}$. By (5.8) and the fact that $\tilde{y}$ satisfies (5.1), we have

$$
\begin{aligned}
(5.14) \quad & K_{M,j_{M+1}}(y_{j_{M+1}}) \\
&= w_{j_{M+1}} \\
&= \sum_{k \in \mathcal{C}_{(j_{M+1})}} \rho_{k,j_{M+1}} \left[ H_k(\tilde{y}_{j_{M+1}}) - H_k\left( \min_{i \in \mathcal{S}(k|j_{M+1})} \tilde{y}_i \right) \right]^+ \\
&\geq \sum_{k \in \mathcal{K}_M^{(j_1, \ldots, j_M)}(j_{M+1})} \rho_{k,j_{M+1}} \left[ H_k(\tilde{y}_{j_{M+1}}) - H_k\left( \min_{i \in \mathcal{S}(k|j_{M+1})} \tilde{y}_i \right) \right]^+ .
\end{aligned}
$$

For $k \in \mathcal{K}_M^{(j_1, \ldots, j_M)}(j_{M+1})$, $\mathcal{S}(k|j_{M+1}) \subseteq \{j_1, \ldots, j_M\}$, so by the induction hypothesis,

$$(5.15) \qquad\qquad \min_{i \in \mathcal{S}(k|j_{M+1})} y_i \leq \min_{i \in \mathcal{S}(k|j_{M+1})} \tilde{y}_i.$$

This, together with (5.14) and the monotonicity of $H_k$, yields

$$
\begin{aligned}
K_{M,j_{M+1}}(y_{j_{M+1}}) &\geq \sum_{k \in \mathcal{K}_M^{(j_1, \ldots, j_M)}(j_{M+1})} \rho_{k,j_{M+1}} \left[ H_k(\tilde{y}_{j_{M+1}}) - H_k\left( \min_{i \in \mathcal{S}(k|j_{M+1})} y_i \right) \right]^+ \\
&= K_{M,j_{M+1}}(\tilde{y}_{j_{M+1}}).
\end{aligned}
$$

Thus, by the monotonicity property of $K_{M,j_{M+1}}$, either

$$(5.16) \qquad\qquad y_{j_{M+1}} \leq \tilde{y}_{j_{M+1}}$$

or

$$(5.17) \qquad \max_{k \in \mathcal{K}_M^{(j_1, \ldots, j_M)}(j_{M+1})} \left( y_k^* \wedge \min_{i \in \mathcal{S}(k|j_{M+1})} y_i \right) \leq \tilde{y}_{j_{M+1}} < y_{j_{M+1}}.$$

However, (5.17) contradicts the definition of $y_{j_{M+1}}$:

$$
y_{j_{M+1}} = K_{M,j_{M+1}}^{-1}(w_{j_{M+1}}) \in \left( -\infty, \max_{k \in \mathcal{K}_M^{(j_1, \ldots, j_M)}(j_{M+1})} \left( y_k^* \wedge \min_{i \in \mathcal{S}(k|j_{M+1})} y_i \right) \right],
$$

so (5.16) holds.  □

LEMMA 5.3.   *The mapping* $\Phi : D \to [0,\infty)^J$ *is one-to-one.*

PROOF.   Let $w = (w_1,\ldots,w_J) \in [0,\infty)^J$ be given and let $\tilde{y} = (\tilde{y}_1,\ldots,\tilde{y}_J) \in D$ be a solution to (5.1). Let $\tilde{\pi} = (\tilde{j}_1,\ldots,\tilde{j}_J) \in \Pi$ be such that $\tilde{y} \in D^{\tilde{\pi}}$, in particular,

$$(5.18) \qquad \tilde{y}_{\tilde{j}_1} \geq \tilde{y}_{\tilde{j}_2} \geq \cdots \geq \tilde{y}_{\tilde{j}_J}.$$

In light of Lemma 5.2, it suffices to show that $\tilde{y}$ and $\tilde{\pi}$ can be constructed as the output $y = (y_1,\ldots,y_J)$, $\pi = (j_1,\ldots,j_J)$, of the algorithm described in Lemma 5.1. Once again, we proceed by induction. We assume that for some $M < J$ and all $m \leq M$ we have chosen in the above-mentioned algorithm $j_m = \tilde{j}_m$ and $y_{j_m} = \tilde{y}_{j_m}$ (for $M = 0$, nothing is assumed). We want to show that it is possible to choose in this algorithm $y_{M+1}$ and $j_{M+1}$ as $\tilde{y}_{M+1}$ and $\tilde{j}_{M+1}$, respectively. By the induction hypothesis, $\mathcal{K}_M^{(j_1,\ldots,j_M)}(j) = \mathcal{K}_M^{\tilde{\pi}}(j)$ for all $j$ and $\mathcal{J}_M^{(j_1,\ldots,j_M)} = \mathcal{J}_M^{\tilde{\pi}}$. In particular, $\tilde{j}_{M+1} \in \mathcal{J}_M^{(j_1,\ldots,j_M)}$ because $\tilde{\pi} \in \Pi$. By assumption,

$$(5.19) \quad w_{\tilde{j}_{M+1}} = \sum_{k \in \mathcal{C}(\tilde{j}_{M+1})} \rho_{k,\tilde{j}_{M+1}} \left[ H_k(\tilde{y}_{\tilde{j}_{M+1}}) - H_k\left( \min_{i \in \mathcal{S}(k|\tilde{j}_{M+1})} \tilde{y}_i \right) \right]^+.$$

If $k \in \mathcal{C}(\tilde{j}_{M+1}) \setminus \mathcal{K}_M^{(j_1,\ldots,j_M)}(\tilde{j}_{M+1})$, then $\mathcal{S}(k|\tilde{j}_{M+1}) \nsubseteq \{\tilde{j}_1,\ldots,\tilde{j}_M\}$ and, hence, by (5.18), $\min_{i \in \mathcal{S}(k|\tilde{j}_{M+1})} \tilde{y}_i \leq \tilde{y}_{\tilde{j}_{M+1}}$. Thus, (5.19) reduces to

$$
\begin{aligned}
w_{\tilde{j}_{M+1}} &= \sum_{k \in \mathcal{K}_M^{(j_1,\ldots,j_M)}(\tilde{j}_{M+1})} \rho_{k,\tilde{j}_{M+1}} \left[ H_k(\tilde{y}_{\tilde{j}_{M+1}}) - H_k\left( \min_{i \in \mathcal{S}(k|\tilde{j}_{M+1})} \tilde{y}_i \right) \right]^+ \\
&= \sum_{k \in \mathcal{K}_M^{(j_1,\ldots,j_M)}(\tilde{j}_{M+1})} \rho_{k,\tilde{j}_{M+1}} \left[ H_k(\tilde{y}_{\tilde{j}_{M+1}}) - H_k\left( \min_{i \in \mathcal{S}(k|\tilde{j}_{M+1})} y_i \right) \right]^+ \\
&= K_{M,\tilde{j}_{M+1}}(\tilde{y}_{\tilde{j}_{M+1}}).
\end{aligned}
$$
$$(5.20)$$

The second equation follows from the fact that

$$(5.21) \quad \mathcal{S}(k|\tilde{j}_{M+1}) \subseteq \{j_1,\ldots,j_M\} = \{\tilde{j}_1,\ldots,\tilde{j}_M\} \qquad \text{for } k \in \mathcal{K}_M^{(j_1,\ldots,j_M)}(\tilde{j}_{M+1})$$

and, hence, by the induction hypothesis, $y_i = \tilde{y}_i$ for $i \in \mathcal{S}(k|\tilde{j}_{M+1})$. But $\tilde{y} \in D^{\tilde{\pi}}$, so

$$\tilde{y}_{\tilde{j}_{M+1}} \leq \max_{k \in \mathcal{K}_M^{(j_1,\ldots,j_M)}(\tilde{j}_{M+1})} y_k^*.$$

By (5.18), (5.21) and the induction hypothesis, for $k \in \mathcal{K}_M^{(j_1,\ldots,j_M)}(\tilde{j}_{M+1})$, we have

$$\tilde{y}_{\tilde{j}_{M+1}} \leq \tilde{y}_{\tilde{j}_M} = \min\{\tilde{y}_{\tilde{j}_1},\ldots,\tilde{y}_{\tilde{j}_M}\} = \min\{y_{j_1},\ldots,y_{j_M}\} \leq \min_{i \in \mathcal{S}(k|\tilde{j}_{M+1})} y_i.$$

Therefore,

$$(5.22) \qquad \tilde{y}_{\tilde{j}_{M+1}} \leq \max_{k \in \mathcal{K}_M^{(j_1, \ldots, j_M)}(\tilde{j}_{M+1})} \left( y_k^* \wedge \min_{i \in \mathcal{S}(k|\tilde{j}_{M+1})} y_i \right).$$

By (5.20) and (5.22),

$$(5.23) \qquad \tilde{y}_{\tilde{j}_{M+1}} = K_{M,\tilde{j}_{M+1}}^{-1}(w_{\tilde{j}_{M+1}}).$$

Since $\tilde{y}$ is a solution to (5.1), for $j \in \mathcal{J}_M^{(j_1, \ldots, j_M)}$, we have

$$
\begin{aligned}
(5.24) \qquad w_j &= \sum_{k \in \mathcal{C}(j)} \rho_{k,j} \left[ H_k(\tilde{y}_j) - H_k \left( \min_{i \in \mathcal{S}(k|j)} \tilde{y}_i \right) \right]^+ \\
&\geq \sum_{k \in \mathcal{K}_M^{(j_1, \ldots, j_M)}(j)} \rho_{k,j} \left[ H_k(\tilde{y}_j) - H_k \left( \min_{i \in \mathcal{S}(k|j)} \tilde{y}_i \right) \right]^+.
\end{aligned}
$$

For $k \in \mathcal{K}_M^{(j_1, \ldots, j_M)}(j)$, $\mathcal{S}(k|j) \subseteq \{j_1, \ldots, j_M\}$, so $y_i = \tilde{y}_i$ for $i \in \mathcal{S}(k|j)$ by the induction hypothesis. Thus, by (5.24),

$$
\begin{aligned}
(5.25) \qquad w_j &\geq \sum_{k \in \mathcal{K}_M^{(j_1, \ldots, j_M)}(j)} \rho_{k,j} \left[ H_k(\tilde{y}_j) - H_k \left( \min_{i \in \mathcal{S}(k|j)} y_i \right) \right]^+ \\
&\geq \sum_{k \in \mathcal{K}_M^{(j_1, \ldots, j_M)}(j)} \rho_{k,j} \left[ H_k(\tilde{y}_{\tilde{j}_{M+1}}) - H_k \left( \min_{i \in \mathcal{S}(k|j)} y_i \right) \right]^+ \\
&= K_{M,j}(\tilde{y}_{\tilde{j}_{M+1}}).
\end{aligned}
$$

The second inequality follows from the fact that $j \in \mathcal{J}_M^{(j_1, \ldots, j_M)}$, so $j \notin \{j_1, \ldots, j_M\}$ and, hence, by (5.18) and the induction hypothesis, $\tilde{y}_j \leq \tilde{y}_{\tilde{j}_{M+1}}$. Relation (5.25) implies

$$(5.26) \quad K_{M,j}^{-1}(w_j) \leq K_{M,j}^{-1}(K_{M,j}(\tilde{y}_{\tilde{j}_{M+1}})) \leq \tilde{y}_{\tilde{j}_{M+1}}, \qquad j \in \mathcal{J}_M^{(j_1, \ldots, j_M)}.$$

By (5.23) and (5.26), $j_{M+1}$ and $y_{M+1}$ in the algorithm of Lemma 5.1 may be chosen as $\tilde{j}_{M+1}$ and $\tilde{y}_{\tilde{j}_{M+1}}$. This is what we wanted to show. $\square$

LEMMA 5.4.   $\lim_{y \in D, \|y\| \to +\infty} \|\Phi(y)\| = +\infty.$

PROOF.   We argue by contradiction. Suppose that the lemma is false. Then there exists a sequence $y_n = (y_1^n, \ldots, y_J^n) \in D$, with $\|y_n\| \to +\infty$ as $n \to +\infty$, and a finite constant $M$ such that

$$(5.27) \qquad \|\Phi(y_n)\| \leq M, \qquad n = 1, 2, \ldots.$$

By taking a subsequence (also denoted by $y_n$), we may assume that for some $\pi = (j_1, \ldots, j_J) \in \Pi$, we have $y_n \in D^\pi$, $n = 1, 2, \ldots$. Let

$$m_0 = \min\{m \in \{1, \ldots, J\} : \{y_{j_m}^n\}_{n=1,2,\ldots} \text{ is an unbounded sequence}\}.$$

By definition the set $D^\pi$ is bounded above in each coordinate. Again extracting a subsequence (still called $y_n$), if necessary, we may assume

$$(5.28) \qquad \lim_{n \to \infty} y_{j_{m_0}}^n = -\infty.$$

Let $w_n = (w_1^n, \ldots, w_J^n) = \Phi(y_n)$, $n = 1, 2, \ldots$. By (5.6) with $m = m_0$, we have

$$(5.29) \quad w_{j_{m_0}}^n = \sum_{k \in \mathcal{K}_{m_0-1}^{(j_1,\ldots,j_{m_0-1})}(j_{m_0})} \rho_{k,j_{m_0}} \left[ H_k\left( y_{j_{m_0}}^n \right) - H_k\left( \min_{i \in \mathcal{S}(k|j_{m_0})} y_i^n \right) \right]^+.$$

[Recall that, by the proof of Lemma 5.3, the permutation constructed in the algorithm of Lemma 5.1 with input $w_n = \Phi(y_n)$, $y_n \in D^\pi$, may be chosen to be $\pi$.] Observe that

$$(5.30) \quad \mathcal{S}(k|j_{m_0}) \subseteq \{j_1, \ldots, j_{m_0-1}\} \qquad \text{for } k \in \mathcal{K}_{m_0-1}^{(j_1,\ldots,j_{m_0-1})}(j_{m_0}),$$

and by the definition of $m_0$,

$$(5.31) \qquad \limsup_{n \to \infty} |y_{j_i}^n| < +\infty, \qquad i = 1, \ldots, m_0 - 1.$$

Relations (5.28)–(5.31) yield $\lim_{n \to \infty} w_{j_{m_0}}^n = +\infty$, because $\lim_{x \to -\infty} H_k(x) = +\infty$ for all $k$. This contradicts (5.27). □

PROPOSITION 5.5. *The mapping* $\Phi : D \to [0, \infty)^J$ *is a homeomorphism of* $D$ *onto* $[0, \infty)^J$.

PROOF. By Lemmas 5.1 and 5.3, it suffices to prove that $\Phi$ is open. The main idea of the proof is to use the one-point (Alexandroff) compactification of $D$ and $[0, \infty)^J$ (see, e.g., [3], pages 92 and 93). Recall that the topology on the one point compactification $\overline{X} = X \cup \{\infty\}$ of a locally compact Hausdorff space $X$ consists of open subsets of $X$ and the complements, in $\overline{X}$, of compact subsets of $X$. Let $\infty$ be a single point not belonging to $R^J$. Let $\overline{D} = D \cup \{\infty\}$ and $\overline{\mathbb{R}}_+^J = [0, \infty)^J \cup \{\infty\}$ be the one-point compactifications of $D$ and $[0, \infty)^J$, respectively. Define $\overline{\Phi} : \overline{D} \to \overline{\mathbb{R}}_+^J$ by

$$\overline{\Phi}(y) \triangleq \begin{cases} \Phi(y), & \text{if } y \in D, \\ \infty, & \text{if } y = \infty. \end{cases}$$

Lemma 5.4 implies the continuity of $\overline{\Phi}$ at $\infty$. Thus, $\overline{\Phi}$ is a continuous mapping of a compact space into a Hausdorff space and, therefore, by Corollary 2 on page 87 of [3], it is closed. In fact, by Lemmas 5.1 and 5.3, $\overline{\Phi}$ is a homeomorphism of $\overline{D}$ onto $\overline{\mathbb{R}}_+^J$. To conclude, let $U$ be an open subset of $D$ and, hence, of $\overline{D}$. Therefore, $\Phi(U) = \overline{\Phi}(U)$ is open in $\overline{\mathbb{R}}_+^J$. But $\Phi(U) \subseteq [0, \infty)^J$, so $\Phi(U)$ is open in $[0, \infty)^J$.

$\square$

## 6. Proofs of the main results.

PROOF OF THEOREM 3.2. By Proposition 5.5, only the second part of Theorem 3.2 needs to be shown. Let us observe that for $j = 1, \ldots, J$, we have [using the convention $(a, b] = \varnothing$ if $a \geq b$]

$$
\begin{aligned}
\widehat{W}_j^{(n)}(t) &= \widehat{\mathcal{W}}_j^{(n)}(t)(\widehat{F}_j^{(n)}(t), \infty) + o(1) \\
&= \sum_{k \in \mathcal{C}(j)} \left[ \widehat{\mathcal{W}}_{k,j}^{(n)}(t)\left( \widehat{F}_j^{(n)}(t), \min_{i \in \mathcal{S}(k|j)} \widehat{F}_i^{(n)}(t) \right] \right. \\
&\qquad\qquad \left. + \widehat{\mathcal{W}}_{k,j}^{(n)}(t)\left( \widehat{F}_j^{(n)}(t) \vee \min_{i \in \mathcal{S}(k|j)} \widehat{F}_i^{(n)}(t), \infty \right) \right] + o(1) \\
&= \sum_{k \in \mathcal{C}(j)} \widehat{\mathcal{W}}_{k,j}^{(n)}(t)\left( \widehat{F}_j^{(n)}(t), \min_{i \in \mathcal{S}(k|j)} \widehat{F}_i^{(n)}(t) \right] + o(1) \\
&= \sum_{k \in \mathcal{C}(j)} \left[ \widehat{\mathcal{V}}_{0,k,j}^{(n)}(t)\left( \widehat{F}_j^{(n)}(t), \min_{i \in \mathcal{S}(k|j)} \widehat{F}_i^{(n)}(t) \right] \right. \\
&\qquad\qquad \left. - \sum_{j_0 \in \mathcal{S}(k|j)} \widehat{\mathcal{W}}_{k,j}^{j_0,(n)}(t)\left( \widehat{F}_j^{(n)}(t), \min_{i \in \mathcal{S}(k|j)} \widehat{F}_i^{(n)}(t) \right] \right] + o(1) \\
&= \sum_{k \in \mathcal{C}(j)} \widehat{\mathcal{V}}_{0,k,j}^{(n)}(t)\left( \widehat{F}_j^{(n)}(t), \min_{i \in \mathcal{S}(k|j)} \widehat{F}_i^{(n)}(t) \right] + o(1) \\
&= \sum_{k \in \mathcal{C}(j)} \rho_{k,j} \left[ H_k(\widehat{F}_j^{(n)}(t)) - H_k\left( \min_{i \in \mathcal{S}(k|j)} \widehat{F}_i^{(n)}(t) \right) \right]^+ + o(1) \\
&= \Phi_j(\widehat{F}_1^{(n)}(t), \ldots, \widehat{F}_J^{(n)}(t)) + o(1).
\end{aligned}
$$

(6.1)

Indeed, the first equality in (6.1) holds by Corollary 4.7. The third one follows from the fact that class $k$ customers with lead times at time $nt$ greater than $\min_{i \in \mathcal{S}(k|j)} F_i^{(n)}(nt)$ have not yet been in service at one of the stations $i \in \mathcal{S}(k|j)$ and, thus, have not yet arrived at station $j$. Similarly, no class $k$ customer with lead time at time $nt$ greater than $F_j^{(n)}(nt)$ has ever been in service at station $j$, so all such customers must be either in queue at station $j$ or at an upstream station $j_0 \in \mathcal{S}(k|j)$. This explains the fourth equality in (6.1). The fifth one follows from the fact that, for every $k \in \mathcal{C}(j)$ and $j_0 \in \mathcal{S}(k|j)$, we have

$$
0 \leq \widehat{\mathcal{W}}_{k,j}^{j_0,(n)}(t)\left( \widehat{F}_j^{(n)}(t), \min_{i \in \mathcal{S}(k|j)} \widehat{F}_i^{(n)}(t) \right]
$$

$$\leq \widehat{\mathcal{W}}_j^{j_0,(n)}(t)(-\infty, \widehat{F}_{j_0}^{(n)}(t)] \Rightarrow 0$$

by Corollary 4.7. Finally, the sixth equation in (6.1) follows from Proposition 4.1 and Lemma 4.6, together with the fact that, by definition, for every station $i \in \mathcal{P}(k)$, $\widehat{F}_i^{(n)}(t) + \sqrt{n}t \geq y_k^*$ and, hence,

$$H_k(\widehat{F}_j^{(n)}(t) + \sqrt{n}t) = H_k\left(\min_{i \in \mathcal{S}(k|j)} \widehat{F}_i^{(n)}(t) + \sqrt{n}t\right) = 0.$$

By Lemma 3.1, $(\widehat{F}_1^{(n)}(t), \ldots, \widehat{F}_J^{(n)}(t)) \in D$ and, by Proposition 5.5, $\Phi$ is a homeomorphism of $D$ onto $[0, \infty)^J$. Thus, (2.25), (3.21) and (6.1), together with the continuous mapping theorem, yield

$$\begin{aligned}
(\widehat{F}_1^{(n)}(t), \ldots, \widehat{F}_J^{(n)}(t)) &= \Phi^{-1}((\widehat{W}_1^{(n)}(t), \ldots, \widehat{W}_J^{(n)}(t)) + o(1)) \\
&\Rightarrow \Phi^{-1}(W_1^*(t), \ldots, W_J^*(t)) = (F_1^*(t), \ldots, F_J^*(t)). \quad \square
\end{aligned}$$

PROPOSITION 6.1. *Let $j = 1, \ldots, J$ and $T > 0$ be given. As $n \to \infty$, both*

$$\sup_{y \in \mathbb{R}} \sup_{0 \leq t \leq T} \left| \widehat{\mathcal{W}}_j^{(n)}(t)(y, \infty) - \sum_{k \in \mathcal{C}(j)} \rho_{k,j} \left[ H_k(y \vee \widehat{F}_j^{(n)}(t)) - H_k\left(\min_{i \in \mathcal{S}(k|j)} \widehat{F}_i^{(n)}(t)\right) \right]^+ \right|$$

*and*

$$\sup_{y \in \mathbb{R}} \sup_{0 \leq t \leq T} \left| \widehat{\mathcal{Q}}_j^{(n)}(t)(y, \infty) - \sum_{k \in \mathcal{C}(j)} \lambda_k \left[ H_k(y \vee \widehat{F}_j^{(n)}(t)) - H_k\left(\min_{i \in \mathcal{S}(k|j)} \widehat{F}_i^{(n)}(t)\right) \right]^+ \right|$$

*converge to zero in probability.*

PROOF. By an argument similar to that used to derive (6.1), we have, uniformly in $0 \leq t \leq T$,

$$\begin{aligned}
\widehat{\mathcal{W}}_j^{(n)}&(t)(y, \infty) \\
&= \widehat{\mathcal{W}}_j^{(n)}(t)(y \vee \widehat{F}_j^{(n)}(t), \infty) + o(1) \\
&= \sum_{k \in \mathcal{C}(j)} \widehat{\mathcal{W}}_{k,j}^{(n)}(t)\left(y \vee \widehat{F}_j^{(n)}(t), \min_{i \in \mathcal{S}(k|j)} \widehat{F}_i^{(n)}(t)\right] + o(1) \\
&= \sum_{k \in \mathcal{C}(j)} \left[ \widehat{\mathcal{V}}_{0,k,j}^{(n)}(t)\left(y \vee \widehat{F}_j^{(n)}(t), \min_{i \in \mathcal{S}(k|j)} \widehat{F}_i^{(n)}(t)\right] \right. \\
&\qquad \left. - \sum_{j_0 \in \mathcal{S}(k|j)} \widehat{\mathcal{W}}_{k,j}^{j_0,(n)}(t)\left(y \vee \widehat{F}_j^{(n)}(t), \min_{i \in \mathcal{S}(k|j)} \widehat{F}_i^{(n)}(t)\right) \right] + o(1) \\
&= \sum_{k \in \mathcal{C}(j)} \widehat{\mathcal{V}}_{0,k,j}^{(n)}(t)\left(y \vee \widehat{F}_j^{(n)}(t), \min_{i \in \mathcal{S}(k|j)} \widehat{F}_i^{(n)}(t)\right] + o(1) \\
&= \sum_{k \in \mathcal{C}(j)} \rho_{k,j} \left[ H_k(y \vee \widehat{F}_j^{(n)}(t)) - H_k\left(\min_{i \in \mathcal{S}(k|j)} \widehat{F}_i^{(n)}(t)\right) \right]^+ + o(1).
\end{aligned}$$

(6.2)

Moreover, the $o(1)$ terms above may be chosen uniformly in $y \in \mathbb{R}$. This needs a justification only for the last equality in (6.2). For $y > y_k^*$, the $k$th terms in the sums in both the fifth and the sixth line of (6.2) are zero. Proposition 4.1 gives a uniform bound for $y_0 \leq y \leq y_k^*$, with $y_0$ arbitrary but fixed. Finally, the uniform bound can be extended to all $y$ by Lemma 4.6. This proves the first part of Proposition 6.1; the proof of the second part is analogous. $\quad\square$

PROOF OF THEOREM 3.3. Let us define a mapping $\psi \colon \mathbb{R}^J \to \mathcal{M}^J$ by $\psi(x) = (\psi_1(x), \dots, \psi_J(x))$, $x = (x_1, \dots, x_J)$, where for $j = 1, \dots, J$ and $B \in \mathcal{B}(\mathbb{R})$,

$$\psi_j(x)(B) = \sum_{k \in \mathcal{C}(j)} \rho_{k,j} \int_{B \cap (x_j, \min_{i \in \mathcal{S}(k|j)} x_i]} (1 - G_k(\xi)) \, d\xi.$$

Observe that $\psi$ is continuous. Indeed, for $j = 1, \dots, J$ and $x, y \in \mathbb{R}^J$, $x = (x_1, \dots, x_J)$, $y = (y_1, \dots, y_J)$, using the fact that $\sum_{k \in \mathcal{C}(j)} \rho_{k,j} = 1$, we have

$$\sup_{B \in \mathcal{B}(\mathbb{R})} |\psi_j(x)(B) - \psi_j(y)(B)|$$

$$\leq \sum_{k \in \mathcal{C}(j)} \rho_{k,j} \int_{(x_j, \min_{i \in \mathcal{S}(k|j)} x_i] \triangle (y_j, \min_{i \in \mathcal{S}(k|j)} y_i]} (1 - G_k(\xi)) \, d\xi$$

$$\leq 2 \max_{l=1,\dots,J} |x_l - y_l|,$$

where $\triangle$ denotes symmetric difference. Therefore, by Theorem 3.2 and the continuous mapping theorem, we have

$$(6.3) \qquad \psi(\widehat{F}_1^{(n)}(t), \dots, \widehat{F}_J^{(n)}(t)) \Rightarrow \psi(F_1^*(t), \dots, F_J^*(t)).$$

For $j = 1, \dots, J$, $t \geq 0$ and $y \in \mathbb{R}$,

$$\psi_j(F_1^*(t), \dots, F_J^*(t))(y, \infty)$$

$$= \sum_{k \in \mathcal{C}(j)} \rho_{k,j} \int_{(y \vee F_j^*(t), \min_{i \in \mathcal{S}(k|j)} F_i^*(t))} (1 - G_k(\xi)) \, d\xi$$

$$= \sum_{k \in \mathcal{C}(j)} \rho_{k,j} \left[ H_k(y \vee F_j^*(t)) - H_k\left( \min_{i \in \mathcal{S}(k|j)} F_i^*(t) \right) \right]^+$$

$$= \mathcal{W}_j^*(t)(y, \infty).$$

This shows that

$$(6.4) \qquad \psi(F_1^*(t), \dots, F_J^*(t)) = \mathcal{W}^*(t).$$

Proposition 6.1 yields

$$(6.5) \quad \sup_{y \in \mathbb{R}} \sup_{0 \le t \le T} |\widehat{\mathcal{W}}_j^{(n)}(t)(y, \infty) - \psi_j(\widehat{F}_1^{(n)}(t), \ldots, \widehat{F}_J^{(n)}(t))(y, \infty)| \xrightarrow{P} 0$$

for every $j = 1, \ldots, J$ and $T > 0$. Combining (6.3), (6.4) and (6.5), we have $(\widehat{\mathcal{W}}_1^{(n)}, \ldots, \widehat{\mathcal{W}}_J^{(n)}) \Rightarrow \mathcal{W}^*$. The proof of $(\widehat{\mathcal{Q}}_1^{(n)}, \ldots, \widehat{\mathcal{Q}}_J^{(n)}) \Rightarrow \mathcal{Q}^*$ is analogous. $\square$

**7. Simulation.** In this section we use simulation methods to assess the predictive value of the theory developed in the previous sections and to provide a simple illustration of the methodology. In the previous sections we considered a sequence of queueing networks, indexed by $n$, whereas here we want to consider a single queueing network. We imagine that this single system is a member of such a sequence of networks corresponding to a large value of $n$, that is, a system with traffic intensities close to one.

Here we show how the theoretical lead-time profile can be constructed when the system occupancy of the $n$th system is given for an EDF network. Suppressing the time variable $t$, we recall that we denote the queue length of class $k$ at station $j$ by $Q_{k,j}^{(n)}$ and its scaled version by

$$\widehat{Q}_{k,j}^{(n)} = \frac{1}{\sqrt{n}} Q_{k,j}^{(n)}.$$

We also denote the workload at station $j$ by $W_j^{(n)}$ and its scaled version by

$$\widehat{W}_j^{(n)} = \frac{1}{\sqrt{n}} W_j^{(n)}.$$

Recall that, by Assumption 2.1, the scaled workloads $(\widehat{W}_1^{(n)}, \ldots, \widehat{W}_J^{(n)})$ converge weakly to $(W_1^*, \ldots, W_J^*)$. The lead-time measure $\mathcal{Q}_{k,j}^{(n)}(y, \infty)$ represents the number of class $k$ customers at station $j$ with lead-time greater than $y$ and

$$\widehat{\mathcal{Q}}_{k,j}^{(n)}(y, \infty) = \frac{1}{\sqrt{n}} \mathcal{Q}_{k,j}^{(n)}(\sqrt{n}y, \infty).$$

Class $k$ customers arrive with lead-time distribution given by

$$(7.1) \qquad\qquad \mathbb{P}(L_k^{i,(n)} \le \sqrt{n}y) = G_k(y).$$

We define $G_k^{(n)}(y) \triangleq G_k(\frac{y}{\sqrt{n}})$ so that

$$(7.2) \qquad\qquad \mathbb{P}(L_k^{i,(n)} \le y) = G_k^{(n)}(y)$$

is the cumulative distribution function of the lead times of class $k$ customers in the $n$th system. The limits of the lead-time measure processes are in terms of the functions $H_k$:

$$(7.3) \qquad H_k(y) \triangleq \int_y^\infty (1 - G_k(\eta))\, d\eta.$$

We also define the function

$$(7.4) \qquad H_k^{(n)}(y) \triangleq \sqrt{n} H_k\left(\frac{y}{\sqrt{n}}\right) = \int_y^\infty (1 - G_k^{(n)}(\eta))\, d\eta.$$

Recall that the frontier at station $j$ is $F_j^{(n)}$, and the scaled frontier is

$$\widehat{F}_j^{(n)} = \frac{1}{\sqrt{n}} F_j^{(n)}.$$

According to Theorem 3.2,

$$(7.5) \qquad (\widehat{F}_1^{(n)}, \ldots, \widehat{F}_J^{(n)}) \Rightarrow (F_1^*, \ldots, F_J^*) \triangleq \Phi^{-1}(W_1^*, \ldots, W_J^*),$$

where $\Phi = (\Phi_1, \ldots, \Phi_J): \mathbb{R}^J \to [0, \infty)^J$ is defined by

$$(7.6) \qquad \Phi_j(y_1, \ldots, y_J) \triangleq \sum_{k \in \mathcal{C}(j)} \rho_{k,j} \left[ H_k(y_j) - H_k\left( \min_{i \in \mathcal{S}(k|j)} y_i \right) \right]^+,$$
$$j = 1, \ldots, J.$$

We also define $\Phi^{(n)} = (\Phi_1^{(n)}, \ldots, \Phi_J^{(n)}): \mathbb{R}^J \to [0, \infty)^J$ by

$$(7.7) \qquad \Phi_j^{(n)}(y_1, \ldots, y_J) \triangleq \sum_{k \in \mathcal{C}(j)} \rho_{k,j}^{(n)} \left[ H_k^{(n)}(y_j) - H_k^{(n)}\left( \min_{i \in \mathcal{S}(k|j)} y_i \right) \right]^+,$$
$$j = 1, \ldots, J.$$

Suppose that, in addition to (2.17), we have

$$\lambda_k - \lambda_k^{(n)} = O\left(\frac{1}{\sqrt{n}}\right) \quad \text{and} \quad \rho_{k,j} - \rho_{k,j}^{(n)} = O\left(\frac{1}{\sqrt{n}}\right)$$

for $k \in \mathcal{C}(j)$, $j = 1, \ldots, J$. Then, by (7.5) and the continuous mapping theorem, we have, for $j = 1, \ldots, J$,

$$\frac{1}{\sqrt{n}} \Phi_j^{(n)}(F_1^{(n)}, \ldots, F_J^{(n)}) = \sum_{k \in \mathcal{C}(j)} \rho_{k,j}^{(n)} \left[ H_k(\widehat{F}_j) - H_k\left( \min_{i \in \mathcal{S}(k|j)} \widehat{F}_i \right) \right]^+$$

$$= \Phi_j(\widehat{F}_1^{(n)}, \ldots, \widehat{F}_J^{(n)}) + O\left(\frac{1}{\sqrt{n}}\right)$$

$$\Rightarrow \Phi_j(F_1^*, \ldots, F_J^*) = W_j^* \approx \frac{1}{\sqrt{n}} W_j^{(n)}.$$

Therefore, $\Phi_j^{(n)}(F_1^{(n)}, \ldots, F_J^{(n)}) \approx W_j^{(n)}$, $j = 1, \ldots, J$. [The difference between these two quantities is $O(1)$, but it is small relative to the number of customers in the system.] Because $\Phi^{(n)}$ has the same functional form as $\Phi$, the proofs in Section 5 apply to $\Phi^{(n)}$ as well as $\Phi$. In particular, $\Phi^{(n)}$ is a homeomorphism of $D^{(n)} \triangleq \sqrt{n}D$ onto $[0, \infty)^J$ (see Proposition 5.5). Therefore,

$$(7.8) \quad (F_1^{(n)}, \ldots, F_J^{(n)}) \approx (\Phi^{(n)})^{-1}(W_1^{(n)}, \ldots, W_J^{(n)}) \triangleq (\overline{F}_1^{(n)}, \ldots, \overline{F}_J^{(n)}).$$

According to Theorem 3.3, for every $y \in \mathbb{R}$ and $j = 1, \ldots, J$, we have, by (7.4), (7.5) and (7.8),

$$
\begin{aligned}
(7.9) \quad \mathcal{Q}_j^{(n)}(y, \infty) &\approx \sqrt{n}\mathcal{Q}_j^*\left(\frac{y}{\sqrt{n}}, \infty\right) \\
&\approx \sum_{k \in \mathcal{C}(j)} \lambda_k\left[H_k^{(n)}(y \vee F_j^{(n)}) - H_k^{(n)}\left(\min_{i \in \mathcal{S}(k|j)} F_i^{(n)}\right)\right]^+ \\
&\approx \sum_{k \in \mathcal{C}(j)} \lambda_k^{(n)}\left[H_k^{(n)}(y \vee \overline{F}_j^{(n)}) - H_k^{(n)}\left(\min_{i \in \mathcal{S}(k|j)} \overline{F}_i^{(n)}\right)\right]^+.
\end{aligned}
$$

In particular,

$$
\begin{aligned}
(7.10) \quad Q_j^{(n)} &= \sqrt{n}\mathcal{Q}_j^*(\mathbb{R}) \\
&\approx \sum_{k \in \mathcal{C}(j)} \lambda_k^{(n)}\left[H_k^{(n)}(\overline{F}_j^{(n)}) - H_k^{(n)}\left(\min_{i \in \mathcal{S}(k|j)} \overline{F}_i^{(n)}\right)\right]^+, \\
&\hspace{6cm} j = 1, \ldots, J.
\end{aligned}
$$

Equations (7.9) and (7.10) indicate that the lead-time profiles can be approximated by a deterministic function in terms of the parameters of the $n$th system, while the knowledge of the index $n$ is not required. The above approximations can be verified by simulation.

### 7.1. *A two station case.*

We consider a simple network with two stations ($J = 2$) and four customer classes ($K = 4$). Flows 1 and 2 visit both stations but in the opposite order, while flows 3 and 4 visit only one station (Figure 1).

### 7.2. *Constant deadline.*

For illustrative purposes, we present the special case in which customers in class $k$ arrive at the system with constant deadline $y_k^*$, that is, $G_k^{(n)}(y) = \mathbb{I}_{[y_k^*, \infty)}(y)$. We also assume that $y_1^* \geq y_2^* \geq y_3^* \geq y_4^*$. We simplify notation by writing $Q_j$, $F_j$, $D$ and $\lambda_k$ in place of $Q_j^{(n)}$, $\overline{F}_j^{(n)}$, $D^{(n)}$ and $\lambda_k^{(n)}$, respectively. It is easy to see that in the case under consideration the sets $\Pi$ and $D$, defined by (3.11)–(3.13), are equal to $\{(1,2),(2,1)\}$ and $D^{(1,2)} \cup D^{(2,1)}$, respectively, where

$$(7.11) \qquad D^{(1,2)} = \{(y_1, y_2): y_1 \geq y_2, y_1 \leq y_1^*, y_2 \leq y_1^*\},$$

$$(7.12) \qquad D^{(2,1)} = \{(y_1, y_2): y_2 \geq y_1, y_1 \leq y_2^*, y_2 \leq y_2^*\}.$$

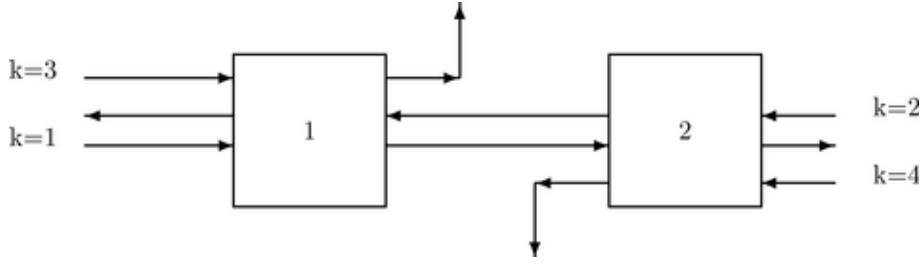

Fig. 1.   *A two node acyclic network with four customer classes.*

Given $Q_1$ and $Q_2$, one can find $F_1$ and $F_2$ by inverting the system of approximate equations (7.10), which in our case reads

$$(7.13) \quad Q_1 \approx \lambda_1(y_1^* - F_1) + \lambda_2[(y_2^* - F_1)^+ - (y_2^* - F_2)^+] + \lambda_3(y_3^* - F_1)^+,$$

$$(7.14) \quad Q_2 \approx \lambda_1(F_1 - F_2)^+ + \lambda_2(y_2^* - F_2)^+ + \lambda_4(y_4^* - F_2)^+.$$

Depending on values of $Q_1$ and $Q_2$, there are five different formulas giving $F_1$ and seven for $F_2$ presented in equations (7.15)–(7.26). These are

$$(7.15) \qquad F_1 \approx y_1^* - \frac{Q_1}{\lambda_1},$$

$$(7.16) \qquad F_1 \approx \frac{\lambda_1 y_1^* + \lambda_2 y_2^* - Q_2 - Q_1}{\lambda_1 + \lambda_2},$$

$$(7.17) \qquad F_1 \approx \frac{\lambda_1 y_1^* + \lambda_3 y_3^* - Q_1}{\lambda_1 + \lambda_3},$$

$$(7.18) \qquad F_1 \approx \frac{\lambda_1 y_1^* + \lambda_2 y_2^* + \lambda_3 y_3^* - Q_1 - Q_2}{\lambda_1 + \lambda_2 + \lambda_3},$$

$$(7.19) \qquad F_1 \approx \frac{\lambda_1 y_1^* + \lambda_3 y_3^* - Q_1}{\lambda_1 + \lambda_2 + \lambda_3} + \frac{\lambda_2(\lambda_2 y_2^* + \lambda_4 y_4^* - Q_2)}{(\lambda_1 + \lambda_2 + \lambda_3)(\lambda_2 + \lambda_4)},$$

$$(7.20) \qquad F_2 \approx y_2^* - \frac{Q_2}{\lambda_2},$$

$$(7.21) \qquad F_2 \approx \frac{\lambda_1 y_1^* - Q_2 - Q_1}{\lambda_1},$$

$$(7.22) \qquad F_2 \approx \frac{\lambda_2 y_2^* + \lambda_4 y_4^* - Q_2}{\lambda_2 + \lambda_4},$$

$$(7.23) \qquad F_2 \approx \frac{\lambda_1 y_1^* + \lambda_2 y_2^* - Q_2 - Q_1}{\lambda_1 + \lambda_2},$$

$$(7.24) \qquad F_2 \approx \frac{\lambda_1 y_1^* + \lambda_2 y_2^* + \lambda_4 y_4^* - Q_2 - Q_1}{\lambda_1 + \lambda_2 + \lambda_4},$$

$$(7.25) \qquad F_2 \approx \frac{\lambda_1 (\lambda_1 y_1^* + \lambda_3 y_3^* - Q_1)}{(\lambda_1 + \lambda_2)(\lambda_1 + \lambda_3)} + \frac{\lambda_2 y_2^* - Q_2}{\lambda_1 + \lambda_2},$$

$$(7.26) \qquad F_2 \approx \frac{\lambda_1 (\lambda_1 y_1^* + \lambda_3 y_3^* - Q_1)}{(\lambda_1 + \lambda_3)(\lambda_1 + \lambda_2 + \lambda_4)} + \frac{\lambda_2 y_2^* + \lambda_4 y_4^* - Q_2}{\lambda_1 + \lambda_2 + \lambda_4}.$$

To describe the function mapping the queue lengths $(Q_1, Q_2)$ to the point $(F_1, F_2) \in D$ satisfying (7.13) and (7.14), we have divided the quadrant $[0, \infty)^2$ in the $(Q_1, Q_2)$-plane into eight regions, I–VIII, which can be seen on Figure 2. Each of these eight regions is mapped onto the corresponding region I'–VIII' of $D$ (i.e., I' is the image of I, etc.) plotted on Figure 3. The regions in $[0, \infty)^2$ are defined by the eight vertices $A$–$H$ shown in Figure 2 and the images of these vertices in $D$ are denoted $A'$–$H'$ in Figure 3.

Table 1 gives the appropriate pair of formulas for the different ranges of values of $(Q_1, Q_2)$, depending on the region in $[0, \infty)^2$ in which the point $(Q_1, Q_2)$ is located.

In the various simulation experiments, we simulate the two node queueing network, as shown in Figure 1. The external inter-arrival times and service times are all assumed to follow exponential distributions. For each simulation run, a particular queue length combination at each node, $Q_i = (Q_{1i}, Q_{2i}, Q_{3i}, Q_{4i})$, $i = 1, 2$, is chosen. Note $Q_{14} = Q_{23} = 0$. The simulation run is initiated with empty queues. Local time is accumulated when the queue length levels at node 1 and node 2 are exactly equal to $Q_1$ and $Q_2$, respectively. At the instant the local time reaches a pre-specified value, 10 for the results presented in this paper, the lead-time profiles for all the customer at each node are recorded and the local time counter is reset to zero. The simulation continues until 50 lead-time profiles at each node are recorded.

| Region | Formula for $F_1$ | Formula for $F_2$ |
|--------|-------------------|-------------------|
| I | (7.15) | (7.21) |
| II | (7.16) | (7.20) |
| III | (7.15) | (7.23) |
| IV | (7.15) | (7.24) |
| V | (7.18) | (7.20) |
| VI | (7.17) | (7.25) |
| VII | (7.19) | (7.22) |
| VIII | (7.17) | (7.26) |

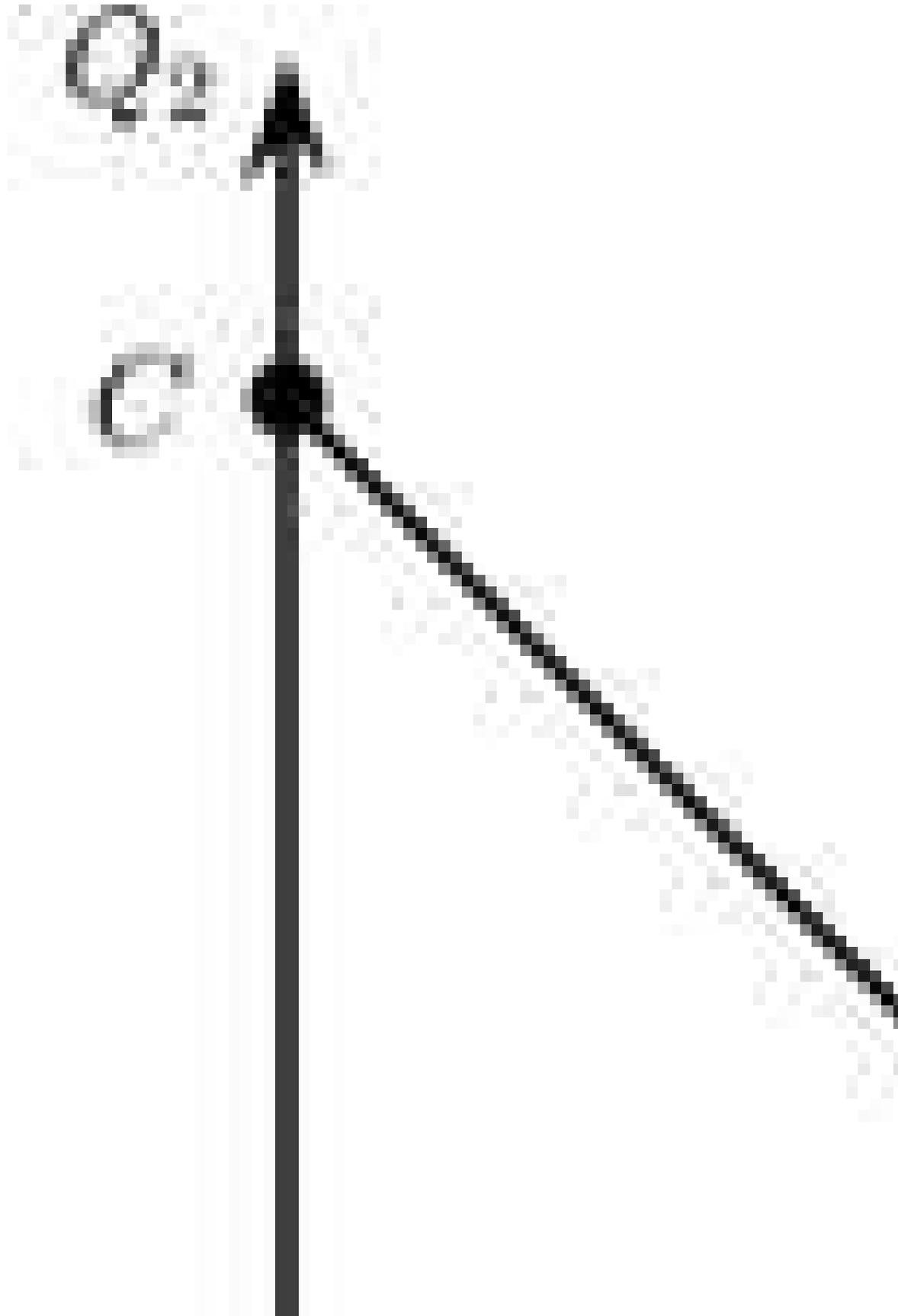

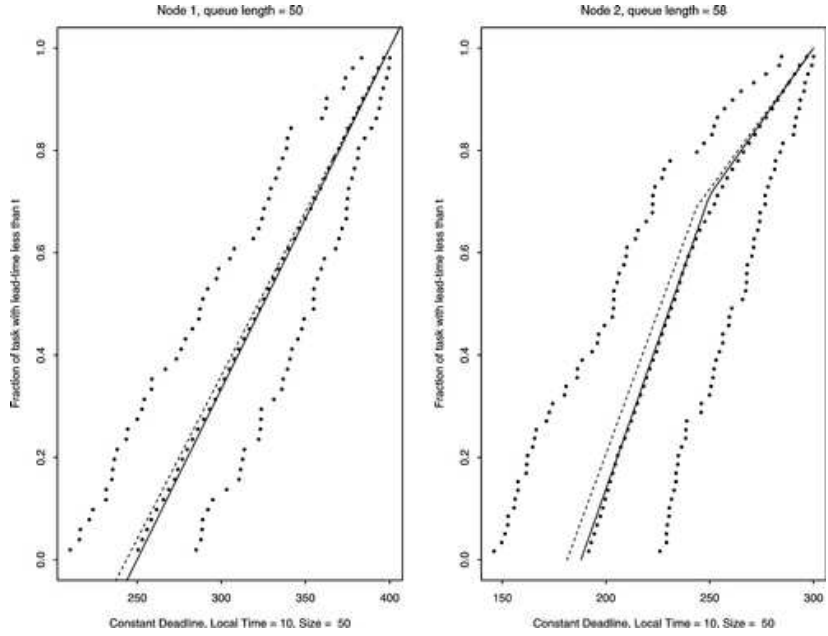

Fig. 4. *Profiles: Mean, Max, Min and Theory,* $Q_1 = (50, 0, 0, 0)$, $Q_2 = (20, 38, 0, 0)$, $D = (400, 300, 200, 100)$.

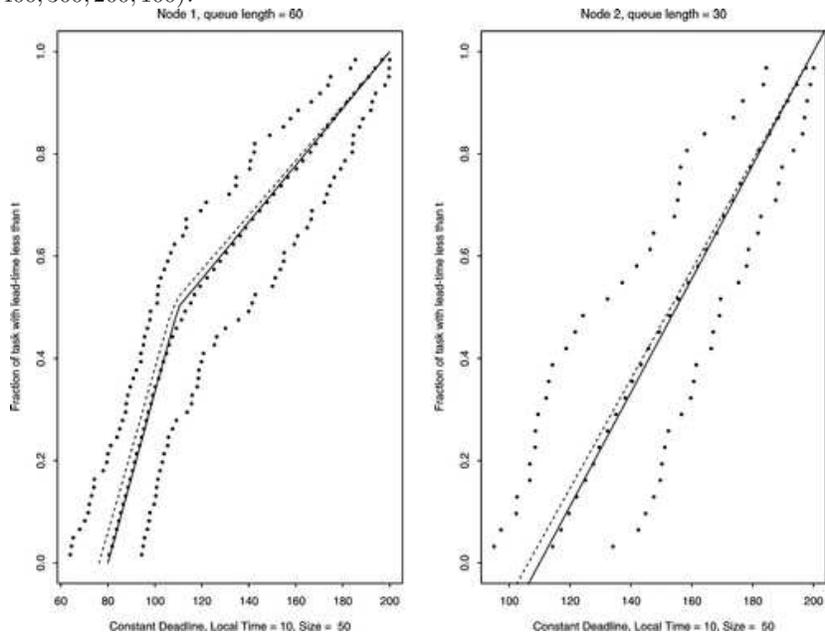

Fig. 5. *Profiles: Mean, Max, Min and Theory,* $Q_1 = (40, 10, 10, 0)$, $Q_2 = (0, 30, 0, 0)$, $D = (200, 200, 110, 100)$.

The empirical lead-time profiles are expressed in the form of empirical lead-time c.d.f.'s.

In the simulation we set $\lambda_1 = \lambda_2 = \lambda_3 = \lambda_4 = 0.32$, and $\mu_{k,j} = 1$ for all $k, j$, so that the total traffic intensity at each node is equal to 0.96. We consider three cases with different combination of queue lengths and end-to-end deadlines. The end-to-end constant deadline for the four customer classes are denoted by $D = (D_1, D_2, D_3, D_4)$. In Figures 4–6, the left-most dots indicate the pointwise minimum empirical cumulative distribution function of the lead-time profile for these 50 samples, the right-most dots indicate the pointwise maximum, and the central dots are the average. As a function of $y$, the ratio $(Q_j^{(n)} - \mathcal{Q}_j^{(n)}(y, \infty))/Q_j^{(n)}$ [where $\mathcal{Q}_j^{(n)}(y, \infty)$ is approximated by the right hand-side of (7.9)] is plotted as a dashed curve in these figures. We obtained the solid curves in Figures 4–6, by replacing $\lambda_j$ in the right-hand sides of (7.9) and (7.10) by $\lambda_j/0.96$. This normalization by the total traffic intensity causes the theory to have better predictive value. Indeed, with this normalization the theoretical cumulative distribution functions and the pointwise average empirical cumulative distribution functions are in almost perfect agreement.

In Figures 4–6, the choice of $Q_1$ and $Q_2$ is made to illustrate different profile compositions at the two nodes. Figure 4 shows the profiles of the case when

$$Q_1 = (50, 0, 0, 0), \qquad Q_2 = (20, 38, 0, 0)$$

and

$$D = (400, 300, 200, 100).$$

In this case only flow 1 is present at node 1, while only flows 1 and 2 are present at node 2. $(F_1, F_2)$ is solved by (7.15) and (7.23) (region III).

In Figure 5 the queue length levels are set at

$$Q_1 = (40, 10, 10, 0), \qquad Q_2 = (0, 30, 0, 0)$$

and the deadlines are

$$D = (200, 200, 110, 100).$$

Here, all three flows are present at node 1, while only flow 2 is present at node 2. In this case $(F_1, F_2)$ is solved by (7.18) and (7.20) (region V).

The final case in Figure 5 shows the profiles when

$$Q_1 = (50, 0, 0, 0), \qquad Q_2 = (50, 0, 0, 0)$$

and

$$D = (500, 100, 100, 100).$$

Only flow 1 is present at both nodes 1 and 2. In this case $(F_1, F_2)$ is solved by (7.15) and (7.21) (region I).

In each of these cases the figures show the excellent predictive accuracy of the theory.

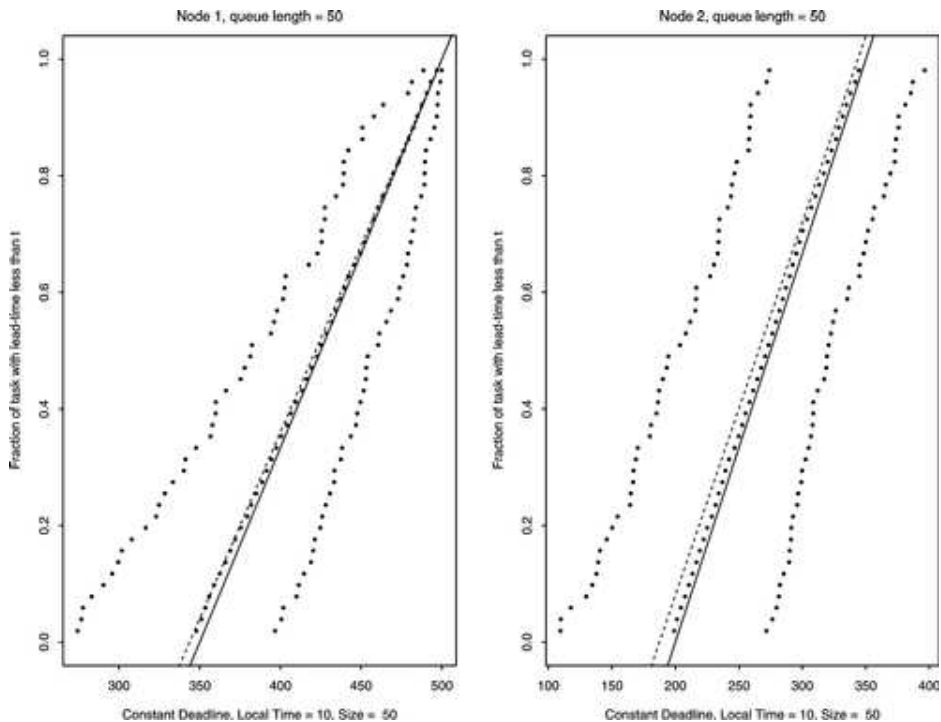

Fig. 6. *Profiles: Mean, Max, Min and Theory,* $Q_1 = (50, 0, 0, 0)$, $Q_2 = (50, 0, 0, 0)$, $D = (500, 100, 100, 100)$.

## REFERENCES

[1] BILLINGSLEY, P. (1986). *Probability and Measure,* 2nd ed. Wiley, New York. MR830424

[2] BILLINGSLEY, P. (1999). *Convergence of Probability Measures,* 2nd ed. Wiley, New York. MR1700749

[3] BOURBAKI, N. (1966). *General Topology* **1**. Addison–Wesley, Reading, MA.

[4] BRAMSON, M. (2004). Stability of earliest-due-date, first-served queueing networks. *Queueing Syst. Theory Appl.* To appear. MR1865459

[5] DOYTCHINOV, B., LEHOCZKY, J. P. and SHREVE, S. E. (2001). Real-time queues in heavy traffic with earliest-deadline-first queue discipline. *Ann. Appl. Probab.* **11** 332–379. MR1843049

[6] ETHIER, S. N. and KURTZ, T. G. (1985). *Markov Processes*: *Characterization and Convergence.* Wiley, New York. MR838085

[7] HARRISON, J. M. (1995). Balanced fluid models of multiclass queueing networks: A heavy traffic conjecture. In *Stochastic Networks* (F. P. Kelly and R. Williams, eds.) 1–20. Springer, New York. MR1381003

[8] KRUK, L., LEHOCZKY, J. P., SHREVE, S. E. and YEUNG, S. N. (2003). Multiple-input heavy-traffic real-time queues. *Ann. Appl. Probab.* **13** 54–99. MR1951994

[9] LEHOCZKY, J. P. (1997). Using real-time queueing theory to control lateness in real-time systems. *Performance Evaluation Review* **25** 158–168.

[10] Lehoczky, J. P. (1998). Real-time queueing theory. In *Proceedings of the IEEE Real-Time Systems Symposium* 186–195.

[11] Lehoczky, J. P. (1998). Scheduling communication networks carrying real-time traffic. In *Proceedings of the IEEE Real-Time Systems Symposium* 470–479.

[12] Markowitz, D. M. and Wein, L. M. (2001). Heavy traffic analysis of dynamic cyclic policies: A unified treatment of the single machine scheduling problem. *Oper. Res.* **49** 246–270. MR1825137

[13] Peterson, W. P. (1991). A heavy traffic limit theorem for networks of queues with different customer types. *Math. Oper. Res.* **16** 90–118. MR1106792

[14] Prokhorov, Yu. (1956). Convergence of random processes and limit theorems in probability theory. *Theory Probab. Appl.* **1** 157–214. MR84896

[15] Van Mieghem, J. A. (1995). Dynamic scheduling with convex delay costs: The generalized $c\mu$ rule. *Ann. Appl. Probab.* **5** 809–833. MR1359830

[16] Williams, R. J. (1998). Diffusion approximations for open multiclass queueing networks: Sufficient conditions involving state space collapse. *Queueing Systems Theory Appl.* **30** 27–88. MR1663759

[17] Yeung, S. N. and Lehoczky, J. P. (2002). Real-time queueing networks in heavy traffic with EDF and FIFO queue discipline. Working paper, Dept. Statistics, Carnegie Mellon Univ.

L. Kruk
Institute of Mathematics
Maria Curie-Sklodowska University
Pl. Marii Curie-Sklodowskiej 1
20-031 Lublin
Poland
e-mail: lkruk@hektor.umcs.lublin.pl

J. Lehoczky
Department of Statistics
Carnegie Mellon University
USA
e-mail: shreve@cmu.edu

S. Shreve
Department of Mathematical Sciences
Carnegie Mellon University
Pittsburgh, Pensylvania 15213-3890
USA
e-mail: jpl@stat.cmu.edu

S.-N. Yeung
AT&T Laboratories
180 Park Avenue
Florham Park, New Jersey
USA
e-mail: syeung@homer.att.com